\documentclass[a4paper]{article}

\pdfoutput=1

\usepackage{epic}
\usepackage{eepic}
\usepackage{amsmath}
\usepackage{amssymb}
\usepackage{theorem}
\usepackage{enumerate}
\usepackage{tikz-cd}
\usepackage{mathtools}

\let\OLDthebibliography\thebibliography
\renewcommand\thebibliography[1]{
	\OLDthebibliography{#1}
	\setlength{\itemsep}{0pt}
}

\newcommand{\al}{\alpha}
\newcommand{\bt}{\beta}

\newcommand{\C}{\mathbb C}
\newcommand{\Q}{\mathbb Q}
\newcommand{\R}{\mathbb R}
\newcommand{\Z}{\mathbb Z}
\newcommand{\la}{\langle}
\newcommand{\ra}{\rangle}

\newcommand{\inv}{\operatorname{inv}}

\newcommand{\sign}{\operatorname{sign}}
\newcommand{\odd}{\operatorname{oddity}}
\newcommand{\tr}{\operatorname{tr}}
\newcommand{\rk}{\operatorname{rk}}
\newcommand{\q}{\operatorname{q}}

\newcommand{\SL}{\operatorname{SL}}
\newcommand{\Mp}{\operatorname{Mp}}
\newcommand{\Orth}{\operatorname{O}}
\newcommand{\Rep}{\operatorname{Re}}
\newcommand{\Imp}{\operatorname{Im}}

\newcommand{\eop}{\hspace*{\fill} $\Box$}

\theoremstyle{break}

\newtheorem{thm}{Theorem}[section]
\newtheorem{prp}[thm]{Proposition}
\newtheorem{cor}[thm]{Corollary}

\begin{document}

\begin{center}
{\Large\bf The invariants of the Weil representation\\[2mm]
  of $\SL_2(\Z)$}\\[10mm]
Manuel K.-H.\ M\"{u}ller, Nils R.\ Scheithauer,\\
Technische Universit\"at Darmstadt, Darmstadt, Deutschland
%Technical University of Darmstadt, Darmstadt, Germany
\end{center}

\vspace*{1.2cm}
\noindent
The transformation behaviour of the vector-valued theta function of a positive-definite even lattice under the metaplectic group $\Mp_2(\Z)$ is described by the Weil representation.
%This representation plays a prominent role in the theory of automorphic forms.
We show that the invariants of this representation are induced from $5$ fundamental invariants.
As an application we
% derive a dimension formula for weight-2 cusp forms for the Weil representation and
give simple generating sets for Jacobi forms of singular weight.

\vspace*{1.2cm}
\begin{tabular}{rl}
1  & Introduction \\
2  & Discriminant forms \\
3  & The Weil representation \\
4  & Discriminant forms of prime level \\
5  & Some $2$-adic exercises \\
6  & Induction \\
7  & The main theorem \\
8  & Applications
\end{tabular}
\vspace*{10mm}

\section{Introduction}

In \cite{W} Weil constructed a representation of the metaplectic group $\Mp_{2n}$, which plays a prominent role in the theory of automorphic forms, called the Weil representation. In the special case of $\Mp_2(\Z)$ it describes the transformation behaviour of the vector-valued theta function of a positive-definite even lattice. For some applications it is important to have an explicit description of the invariants of the Weil representation of $\Mp_2(\Z)$. For example the space of Jacobi forms of lattice index $L$
%where $L$ is a positive definite even lattice
and singular weight is naturally isomorphic to the space of invariants $\C[L'/L]^{\Mp_2(\Z)}$ (cf.\ \cite{Sk2}). If the corresponding discriminant form possesses self-dual isotropic subgroups, the invariants have been described by Skoruppa (cf.\ \cite{Sk2}, \cite{Bi} and \cite{NRS}). They are generated by the characteristic functions of these groups. In the present paper we give a complete description of the invariants for arbitrary discriminant forms. We show that they are induced from $5$ fundamental invariants.

We describe our results in more detail. A discriminant form is a finite abelian group $D$ with a non-degenerate quadratic form $\q : D \to \Q/\Z$. The level of $D$ is the smallest positive integer $N$ such that $N\q(\gamma) = 0 \! \mod 1$ for all $\gamma \in D$. The square class of $D$ is square if $|D|$ is a square and non-square otherwise. Every discriminant form can be realised as the quotient $L'/L$ where $L$ is an even lattice with dual $L'$. The signature of $L$ is unique modulo $8$. We can even assume that $L$ is positive-definite. The vector-valued theta function of $L$ is a function on the upper halfplane ${\cal H}$ with values in the group algebra $\C[D]$ defined by 
\[ \theta(\tau) = 
\sum_{\gamma \in D} \theta_{\gamma}(\tau) e^{\gamma} \]
where $\theta_{\gamma}(\tau) = \sum_{\al \in \gamma+L} q^{\al^2/2}$. Here $q^{\al^2/2} = e^{2\pi i \tau \al^2/2} = e ( \tau \al^2/2)$.
The Poisson summation formula
%The Jacobi transformation formula
implies that $\theta$ transforms as a vector-valued modular form of weight $\rk(L)/2$ under the metaplectic cover $\Mp_2(\Z)$ of $\SL_2(\Z)$. The corresponding representation $\rho_D$ of $\Mp_2(\Z)$ on the group algebra $\C[D]$ is called the Weil representation of $\Mp_2(\Z)$. The non-trivial element in the kernel of the covering map $\Mp_2(\Z) \to \SL_2(\Z)$ acts as $(-1)^{\sign(D)}$ so that $\C[D]^{\Mp_2(\Z)}$ is trivial if $D$ has odd signature. Hence we can restrict to the case that
% the signature of $D$ is even
$D$ has even signature when we study the subspace of invariants. Then the Weil representation $\rho_D$ 
descends to a representation of $\SL_2(\Z)$. 
%factors through $\SL_2(\Z)$. 

Now let $D$ be a discriminant form of even signature and level $N$. Then the Weil representation of $\SL_2(\Z)$ factors through the finite group 
$\SL_2(\Z)/\Gamma(N) \cong \SL_2(\Z/N\Z)$. 
% $\SL_2(\Z)/\Gamma(N) = \SL_2(\Z/N\Z)$.
Hence we can project onto the subspace of invariants by averaging. We define the map
\[  \inv_D : \C[D] \to \C[D]   \]
with
\[  \inv_D(e^{\gamma}) = \frac{1}{|\Gamma(N) \backslash \SL_2(\Z)|}
  \sum_{M \in \Gamma(N) \backslash \SL_2(\Z)} \rho_D(M^{-1}) e^{\gamma} \, . \]
We derive explicit formulas for $\inv_D(e^{\gamma})$ and for the dimension of the space of invariants (see Theorems \ref{invexplicit} and \ref{dimformula}). We evaluate the formulas for discriminant forms of prime level and for some 2-adic discriminant forms of level 4 and 8 (see Sections \ref{wims} and \ref{zwei}).

Let $H$ be an isotropic subgroup of $D$. Then $H^{\perp}/H$ is a discriminant form of the same signature as $D$ and of order $|H^{\perp}/H| = |D|/|H|^2$. There is an isotropic lift $\uparrow_H^D : \C[H^{\perp}/H] \to \C[D]$ which commutes with the corresponding Weil representations (see Section \ref{heidegger}). In particular $\uparrow_H^D$ maps invariants to invariants. We will use these maps to construct all invariants from certain fundamental invariants. 

Let $N=\prod_{p|N} p^{\nu_p}$ be the prime decomposition of $N$. Then $D$ decomposes into the orthogonal sum of
%its $p$-subgroups
$p$-adic discriminant forms
$D = \oplus_{p|N} D_{p^{\nu_p}}$
%\[  D = \bigoplus_{p|N} D_{p^{\nu_p}}   \]
with $D_{p^{\nu_p}} = \{ \gamma \in D \, | \, p^{\nu_p} \gamma = 0 \}$.
We can factorise $\SL_2(\Z/N\Z)$ correspondingly as
$\SL_2(\Z/N\Z) \cong \prod_{p|N} \SL_2(\Z/p^{\nu_p}\Z)$.
%\[  \SL_2(\Z/N\Z) \cong \prod_{p|N} \SL_2(\Z/p^{\nu_p}\Z)  \, .  \]
Then 
\[  \C[D]^{\SL_2(\Z/N\Z)} 
\cong \bigotimes_{p|N} \C[D_{p^{\nu_p}}]^{\SL_2(\Z/p^{\nu_p}\Z)}  \]
so that in order to describe the invariants of the Weil representation of $\SL_2(\Z)$ it suffices to consider $p$-adic discriminant forms.

For this purpose we define $5$ fundamental discriminant forms $D^{x,s}_p$ of square class $x$ and signature $s$. We show that their subspace of invariants is $1$-di\-men\-sio\-nal and determine a generator $i^{x,s}_p$.
% We list them in the following tables.
For odd $p$ the fundamental discriminant forms are given by

\vspace*{1mm}
\[
\renewcommand{\arraystretch}{1.2}
\begin{array}{c|c|c|c}
 D^{x,s}_p  & \text{square class} & \text{signature} & i^{x,s}_p \\[0.5mm] \hline 
   &   &   & \\[-4mm]
 0           & \text{square}           & 0 \, \bmod 8      & e^0 \\
 p^{-4}       & \text{square}           & 4 \, \bmod 8      
                       & (p-1) e^0 - \sum_{\gamma \in M} e^{\gamma} \\
 p^{\epsilon 3} & \text{non-square}       & 0 \, \bmod 2      
                       & \sum_{\gamma \in M^+}e^{\gamma} - \sum_{\gamma \in M^-}e^{\gamma} 
\end{array}
\]

\vspace*{2mm}
\noindent
and for $p=2$ by

\vspace*{-3mm}
\[
\renewcommand{\arraystretch}{1.2}
\begin{array}{c|c|c|c}
 D^{x,s}_p  & \text{square class} & \text{signature} & i^{x,s}_p \\[0.5mm] \hline    &   &   & \\[-4mm]
 0           & \text{square}           & 0 \, \bmod 8      & e^0 \\
 2_{I\!I}^{-4}       & \text{square}           & 4 \, \bmod 8      
                       &  e^0 - \sum_{\gamma \in M} e^{\gamma} \\
 2_t^{+2}4_{I\!I}^{+2}  &  \text{square}  & t = 2 \, \bmod 4 
                       & \sum_{\gamma \in M^+}e^{\gamma} - \sum_{\gamma \in M^-}e^{\gamma} \\
 2_1^{+1} 4_t^{\epsilon} \, 8_{I\!I}^{+2}  &  \text{non-square}  &  1+t = 0 \, \bmod 2 
                       & \sum_{\gamma \in M^+}e^{\gamma} - \sum_{\gamma \in M^-}e^{\gamma} 
\end{array}
\]

\vspace*{2mm}
\noindent
The notation for the discriminant forms will be explained in Section \ref{df}. We wrote $M$ for the set of isotropic elements whose order is equal to the level of $D^{x,s}_p$. In the indicated cases $M$ has a canonical decomposition $M = M^+ \cup M^-$. Our main result is the following (see Theorem \ref{mainth}):

{\em Let $D$ be a discriminant form of even signature $s$, square class $x$ and level $p^l$ where $p$ is a prime. Then the invariants of the Weil representation on $\C[D]$ are generated by the invariants $\uparrow_H^D(i^{x,s}_p)$ where $H$ is an isotropic subgroup of $D$ such that $H^{\perp}/H$ is isomorphic to the discriminant form $D^{x,s}_p$.}

The idea of the proof is to show that for each $\gamma \in D$, $\inv_D(e^{\gamma})$ is a linear combination of isotropic lifts of invariants for suitable isotropic subgroups unless $\inv_D(e^{\gamma})=0$ or $D$ is the fundamental discriminant form $D^{x,s}_p$. Then by induction on the order of $D$ the invariant $\inv_D(e^{\gamma})$ is a linear combination of compositions of isotropic lifts of invariants of the form $i^{x,s}_p$. The statement now follows from the transitivity of the isotropic lift. 

%We remark that
Skoruppa's result corresponds to the case that $D^{x,s}_p$ is trivial.

As an application of our main result we show (see Theorem \ref{Iggypop}):

{\em Let $L$ be a positive-definite even lattice of even rank $n$ and level $N$. For $p|N$ we denote the square class and the signature of the $p$-adic component of $L'/L$ by $x_p$ resp.\ $s_p$. Let $\cal{L}$ be the set of all overlattices $M \supset L$ such that the $p$-adic component of $M'/M$ is isomorphic to $D_p^{x_p,s_p}$ for all $p|N$.
Then the space $J_{n/2,L}$ of Jacobi forms of lattice index $L$ and weight $n/2$ is generated by the functions
\[ \sum_{\gamma \in M'/M} v_{\gamma} \vartheta_{M,\gamma} \]
where $M \in {\cal L}$, $\sum_{\gamma\in M'/M} v_{\gamma} e^\gamma \in \C[M'/M]^{\SL_2(\Z)}$ is the invariant corresponding to the product $\prod_{p|N} i_p^{x_p,s_p}$ and
%$\vartheta_{M, \gamma}(\tau,z) = \sum_{\al \in \gamma + M} e( \tau \al^2/2 + (\al,z) )$
\[  \vartheta_{M, \gamma}(\tau,z) = \sum_{\al \in \gamma + M} e( \tau \al^2/2 + (\al,z) )  \]
is the Jacobi theta function of the coset $\gamma + M$.}

\medskip

The paper is organised as follows. In Section 2 we recall some results about discriminant forms. In Section 3 we recall the Weil representation of $\SL_2(\Z)$ and give an explicit description of the projection map $\inv_D$. We also prove a dimension formula for the subspace of invariants. Then we evaluate these formulas for discriminant forms of prime level. For our main theorem we need some additional results on $2$-adic discriminant forms which we prove in Section 5. Next we recall some properties of the isotropic induction. In Section $7$ we define the fundamental invariants and prove our main theorem. Finally we describe two applications of our results. We determine the dimension of a space of weight-2 cusp forms for the Weil representation and give generators of the space of Jacobi forms of lattice index $L$ and singular weight.

\medskip

We thank P.\ Bieker, M.\ Dittmann, T.\ Driscoll-Spittler, P.\ Kiefer and S.\ Zemel for stimulating discussions and helpful comments. We also thank the referee for carefully reading the manuscript and for suggesting several improvements.

\medskip

Both authors acknowledge support 
% by the LOEWE research unit \emph{Uniformized Structures in Arithmetic and Geometry} and
by the DFG through the CRC \emph{Geometry and Arithmetic of Uniformized Structures}, project number 444845124.

\section{Discriminant forms} \label{df}

In this section we recall some results on discriminant forms (cf.\ 
\cite{AGM}, \cite{B2}, \cite{CS}, \cite{N}, \cite{S2} and \cite{Sk2}). 

\medskip

A discriminant form is a finite abelian group $D$ with a quadratic form $\q : D \to \Q/\Z$ such that $(\bt,\gamma) = \q(\bt + \gamma)- \q(\bt) - \q(\gamma) \! \mod 1$ is a non-degenerate symmetric bilinear form. The level of $D$ is the smallest positive integer $N$ such that $N\q(\gamma) = 0 \! \mod 1$ for all $\gamma \in D$. The square class of $D$ is square if $|D|$ is a square and non-square otherwise. We denote by $\Orth(D)$ the group of automorphisms of $D$ which preserve the quadratic form $\q$. (More generally a homomorphism between discriminant forms is a group homomorphism which preserves the quadratic forms.) 

If $L$ is an even lattice, then $L'/L$ is a discriminant form with the quadratic form given by $\q(\gamma) = \gamma^2/2 \! \mod 1$. Conversely every discriminant form can be obtained in this way. The corresponding lattice can be chosen to be positive-definite. (This can be seen as follows. Choose a lattice realising the given discriminant form. Add a number of $E_8$-lattices and then split off sufficiently many hyperbolic planes $I\!I_{1,1}$.) The signature $\sign(D) \in \Z/8\Z$ of a discriminant form $D$ is defined as the signature modulo $8$ of any even lattice with that discriminant form.

Every discriminant form decomposes into a sum of Jordan components and every Jordan component can be written as a sum of indecomposable Jordan components (usually not uniquely). The possible non-trivial Jordan components are the following:

Let $q>1$ be a power of an odd prime $p$. The non-trivial $p$-adic Jordan components of exponent $q$ are $q^{\pm n}$ for $n\geq 1$. The indecomposable components are $q^{\pm 1}$, generated by an element $\gamma$ with $q\gamma = 0, \, \q(\gamma)= a/q \!\mod 1$ where $a$ is an integer with $\big( \frac{2a}{p} \big) = \pm 1$. These components all have level $q$. The $p$-excess is given by 
$\text{$p$-excess}(q^{\pm n}) = n(q-1) + 4k \mod 8$ where $k=1$, if $q$ is not a square and the exponent is $-n$, and $k=0$ otherwise. We define $\gamma_p(q^{\pm n}) = e(-\text{$p$-excess}(q^{\pm n})/8)$.

Let $q>1$ be a power of $2$. The non-trivial even $2$-adic Jordan components of exponent $q$ are $q^{\pm 2n}=q_{I\!I}^{\pm 2n}$ for $n\geq 1$. The indecomposable components are $q_{I\!I}^{\pm 2}$ generated by two elements $\gamma$ and $\delta$ with $q \gamma = q \delta = 0, \, (\gamma,\delta)= 1/q \!\mod 1$ and $\q(\gamma) = \q(\delta) = 0 \! \mod 1$ for $q_{I\!I}^{+2}$ and $\q(\gamma) = \q(\delta) = 1/q \!\mod 1$ for $q_{I\!I}^{-2}$. These components all have level $q$. The oddity is given by 
$\odd (q_{I\!I}^{\pm 2n}) = 4k \! \mod 8$ 
with $k=1$, if $q$ is not a square and the exponent is $-2n$, and 
%$k=1$, if $\sq(q)=-1$ and the exponent is $-n$, and 
$k=0$ otherwise. We define $\gamma_2(q_{I\!I}^{\pm 2n}) = e(\odd (q_{I\!I}^{\pm 2n})/8)$.

Let $q>1$ be a power of $2$. The non-trivial odd $2$-adic Jordan components of exponent $q$ are $q_{t}^{\pm n}$ with $n\geq 1$ and $t \in \Z/8\Z$. If $n=1$, then $\pm=+$ implies $t=\pm 1 \!\mod 8$ and $\pm=-$ implies $t=\pm 3 \!\mod 8$. If $n=2$, then $\pm=+$ implies $t=0$ or $\pm 2 \!\mod 8$ and $\pm=-$ implies $t=4$ or $\pm 2 \!\mod 8$. For any $n$ we have $t=n \!\mod 2$. The indecomposable components are $q_{t}^{\pm 1}$ where $\big( \frac{t}{2} \big) = \pm 1$ (recall that $\big( \frac{t}{2} \big) = +1$ if $t = \pm 1 \bmod 8$ and $\big( \frac{t}{2} \big) = -1$ if $t = \pm 3 \bmod 8$) generated by an element $\gamma$ with $q\gamma =0, \, \q(\gamma) = t/2q \!\mod 1$. These components all have level $2q$. The oddity is given by 
$ \odd(q_{t}^{\pm n}) = t + 4k \! \mod 8$ 
with $k=1$, if $q$ is not a square and the exponent is $-n$, and 
%$k=1$, if $\sq(q)=-1$ and the exponent is $-n$, and 
$k=0$ otherwise. We define $\gamma_2(q_{t}^{\pm n}) = e(\odd(q_{t}^{\pm n})/8)$.

The sum of two Jordan components with the same prime power $q$ is given by multiplying the signs, adding the ranks and, if any components have a subscript $t$, adding the subscripts $t$. Isomorphic discriminant forms can have different $2$-adic symbols.

Let $D$ be a discriminant form. Then 
\[ \sign(D) + \sum_{p \geq 3} \text{$p$-excess}(D) = \odd(D) \mod 8 \]
respectively
\[ \prod \gamma_p(D) = e(\sign (D)/8)\, .  \]
We will also use 
\[ e(\odd(D)/4) = \left( \frac{-1}{|D|} \right) e(\sign (D)/4) \, . \]

Let $c$ be an integer. Then $c$ acts by multiplication on $D$ and we have an exact sequence
$0 \to D_c \to D \to D^c \to 0$
%\[ 0 \to D_c \to D \to D^c \to 0  \]
where $D_c$ is the kernel and $D^c$ the image of this map. Note that $D^c$ is the orthogonal complement of $D_c$. The set $D^{c*} = \{ \gamma \in D \, | \, c\q(\alpha) + (\alpha,\gamma) = 0 \! \mod 1 \, \text{ for all } \alpha \in D_c \}$ is a coset of $D^c$. After a choice of Jordan decomposition of $D$ there is a canonical coset representative $x_c \in D$ with $2x_c=0$. We can write $\gamma \in D^{c*}$ as $\gamma = x_c + c \mu$. Then $\q_c(\gamma) = c\q(\mu) + x_c \mu \mod 1$ defines a map
$\q_c : D^{c*} \to \Q/\Z$
%\[  \q_c : D^{c*} \to \Q/\Z   \]
which is invariant under $\Orth(D)$ (see Theorem 3.9 in \cite{S2}). 

%We define $D^{c*}$ as the set of elements $\al \in D$ satisfying 
%\[ c\gamma^2/2 + \al \gamma = 0 \mod 1 \]
%for all $\gamma \in D_c$. Then 

\medskip

We describe the number of elements of a given norm in $p$-elementary discriminant forms. To simplify the notation we define for $x \in \Q/\Z$
\[ 
\delta(x) = 
\begin{cases}
        \, 1 & \text{if $x = 0    \! \mod 1$,} \\
        \, 0 & \text{if $x \neq 0 \! \mod 1$.}
\end{cases}
\]
For odd primes we have (see Proposition 3.2 in \cite{S1})

\begin{prp} \label{numberofelementspodd}
Let $p$ be an odd prime. Then the number $N(p^{\epsilon n},j)$ of elements of norm $j/p \! \mod 1$ in the discriminant form $p^{\epsilon n}$ is given by
\[ N(p^{\epsilon n},j) = 
\begin{cases}
\, {\displaystyle p^{n-1} + \epsilon \left( \frac{-1}{p} \right)^{n/2} \big( p \delta(j/p)-1 \big) p^{(n-2)/2} }  
& \text{if $n$ is even}, \\[3mm]
\, {\displaystyle p^{n-1} + \epsilon \left( \frac{-1}{p} \right)^{(n-1)/2}\left( \frac{2}{p} \right) \left( \frac{j}{p} \right) p^{(n-1)/2} }
& \text{if $n$ is odd}.
\end{cases}
\]
\end{prp}

In the level $2$ case we have (see Proposition 3.1 in \cite{S1})

\begin{prp} \label{numberofelements2even}
The number of elements of norm $j/2 \! \mod 1$ in $2_{I\!I}^{\epsilon n}$ is given by
\[   N(2_{I\!I}^{\epsilon n},j) = 2^{n-1} + \epsilon(-1)^j 2^{(n-2)/2}  \, . \]
\end{prp}

Finally for the level $4$ case 

\begin{prp} \label{numberofelements2odd}
The number of elements of norm $j/4 \! \mod 1$ in $2_t^{\epsilon n}$ is given by
\[
N(2_t^{\epsilon n},j) = 
\begin{cases}
\, \displaystyle 2^{n-2} + \epsilon \left( \frac{t}{2} \right) 2^{(n-3)/2}              & \text{if $j = 0 \mod 4$,} \\[3mm]
\, \displaystyle 2^{n-2} - \epsilon \left( \frac{t}{2} \right) 2^{(n-3)/2}              & \text{if $j = 2 \mod 4$,} \\[3mm]
\, \displaystyle 2^{n-2} + \epsilon \left( \frac{t}{2} \right) (-1)^{(t-1)/2}2^{(n-3)/2}  & \text{if $j = 1 \mod 4$,} \\[3mm]
\, \displaystyle 2^{n-2} - \epsilon \left( \frac{t}{2} \right) (-1)^{(t-1)/2}2^{(n-3)/2}  & \text{if $j = 3 \mod 4$}
\end{cases}
\]
if $n$ is odd and by
\[
N(2_t^{\epsilon n},j) = 
\begin{cases}
\, \displaystyle   2^{n-2} + \epsilon \delta(t/4) \left( \frac{t-1}{2} \right) 2^{(n-2)/2}     & \text{if $j = 0 \mod 4$,} \\[3mm]
\, \displaystyle   2^{n-2} - \epsilon \delta(t/4) \left( \frac{t-1}{2} \right) 2^{(n-2)/2}     & \text{if $j = 2 \mod 4$,} \\[3mm]
\, \displaystyle   2^{n-2} + \epsilon \delta((t+2)/4) \left( \frac{t-1}{2} \right) 2^{(n-2)/2} & \text{if $j = 1 \mod 4$,} \\[3mm]
\, \displaystyle   2^{n-2} - \epsilon \delta((t+2)/4) \left( \frac{t-1}{2} \right) 2^{(n-2)/2} & \text{if $j = 3 \mod 4$}
\end{cases}
\]
if $n$ is even.
\end{prp}
{\em Proof:} As in the previous cases this can be proved by induction on $n$. \eop

\section{The Weil representation}

In this section we recall the Weil representation of $\SL_2(\Z)$. Then we describe the projection on the subspace of invariants and derive a formula for the dimension of the space of invariants.

\medskip

Let $D$ be a discriminant form with quadratic form $\q : D \to \Q/\Z$ of even signature. 
%Let $N$ be a positive integer such that the level of $D$ divides $N$. 
We define a scalar product on the group ring $\C [D]$ which is linear in the first and antilinear in the second variable by $(e^{\gamma},e^{\bt})= \delta_{\gamma \bt}$. There is a unitary action $\rho_D$ of the group $\Gamma = \SL_2(\Z)$ on $\C[D]$ given by
\begin{align*} 
\rho_D(T) e^{\gamma}  & = e(-\q(\gamma))\, e^{\gamma} \\
\rho_D(S) e^{\gamma}  & = \frac{e(\sign(D)/8)}{\sqrt{|D|}}
                  \sum_{\bt\in D} e((\gamma,\bt))\, e^{\bt}  
\end{align*}
where 
$S = \left( \begin{smallmatrix} 0 & -1 \\ 1 & 0 \end{smallmatrix} \right)$ and 
$T = \left( \begin{smallmatrix} 1 &  1 \\ 0 & 1 \end{smallmatrix} \right)$
are the standard generators of $\Gamma$. This rep\-re\-sen\-tation is called the Weil representation. With this definition the theta function of a positive-definite even lattice of even rank transforms under the dual Weil representation. (The definition here is dual to the one used in \cite{B1}.)

The element $Z=S^2=-1$ acts as 
\[  \rho_D(Z) e^{\gamma} = e(\sign(D)/4) e^{-\gamma}  \, . \]
For a matrix $M = \left( \begin{smallmatrix} a & b \\ c & d \end{smallmatrix} \right) \in \Gamma$ we have
\[ \rho_D(M) e^{\gamma} = \xi(M) \frac{\sqrt{|D_c|}}{\sqrt{|D|}}  
                    \sum_{\bt \in D^{c*}} 
                    e(-a \q_c(\bt)) \, e(-b(\bt,\gamma)) \, e(-bd \q(\gamma)) \,  
                    e^{d\gamma + \bt}
\]
with $\xi(M) = e(\sign(D)/4) \prod \xi_p(M)$. The local factors $\xi_p(M)$ can be expressed in terms of the Jordan components of $D$ (see Theorem 4.7 in \cite{S2}). 

Let $N$ be a positive integer such that the level of $D$ divides $N$. If $c=0 \!\mod N$, the above formula simplifies to
\[ \rho_D(M) e^{\gamma} = 
%\left( \frac{a}{|D|} \right) e\big( (a-1) \odd(D)/8 \, \big) \, e(-bd \q(\gamma)) \, e^{d\gamma} \, . 
\chi_D(a) \, e(-bd \q(\gamma)) \, e^{d\gamma} 
\]
where 
\[ \chi_D(a) = \left( \frac{a}{|D|} \right) e( (a-1) \odd(D)/8 \, )  \]
is a quadratic Dirichlet character modulo $N$. In particular $\Gamma(N)$ acts trivially. 

\medskip

We denote the set of isotropic elements in $D$ by $I$. Let $v = \sum_{\gamma \in D} v_{\gamma} e^{\gamma}$ be an invariant of $\Gamma$. Then the $T$-invariance implies that $v_{\gamma} = 0$ if $\gamma \notin I$. Hence $\dim \C[D]^{\Gamma} \leq |I|$. We give an exact formula below. 
%The formula for the action of $\Gamma_0(N)$ implies $\chi_D(a)v_{\gamma} =  v_{a\gamma}$ if $(a,N)=1$.

\medskip

We recall some properties of $\Gamma(N)$. The group $\Gamma(N)$ is a normal subgroup of $\Gamma$ and has index 
\[  | \Gamma(N) \backslash \Gamma | = N^3 \prod_{p|N} (1-1/p^2)  \]
in $\Gamma$. The number of classes of cusps is
\[  | \Gamma(N) \backslash P | =
  \begin{cases}
  \, 3                            & \text{if } \, N = 2 \, , \\
  \, (N^2/2) \prod_{p|N} (1-1/p^2) & \text{if } \, N \geq 3 \, 
\end{cases}
\]
where $P = \Q \cup \{ \infty \}$. The (classes of) cusps are parametrised by the elements $(a,c)$ of order $N$ in $(\Z/N\Z)^2$ if $N=2$ and by the pairs $\pm (a,c)$ of elements of order $N$ in $(\Z/N\Z)^2$ if $N\geq 3$ (see Lemma 3.8.4 in \cite{DS}). Let $M=\left( \begin{smallmatrix} a & b \\ c & d \end{smallmatrix} \right) \in \Gamma$. Then the cosets of $\Gamma(N) \backslash \Gamma$ sending $\infty$ to $a/c$ can be represented by $MT^n$ if $N=2$ and by $\pm MT^n$ if $N \geq 3$ where in both cases $n$ ranges over a complete set of residues modulo $N$. 

\medskip

Now we describe the projection on the subspace of invariants. We define the map
\[  \inv_D : \C[D] \to \C[D]   \]
by
\[  \inv_D(e^{\gamma}) = \frac{1}{|\Gamma(N) \backslash \Gamma|}
  \sum_{M \in \Gamma(N) \backslash \Gamma} \rho_D(M^{-1}) e^{\gamma} \, . \]
It maps onto 
%Its image is contained in 
the subspace of invariants $\C[D]^{\Gamma}$. Since $\rho_D$ is unitary, we have
\[  (\inv_D(v),w) = (v, \inv_D(w))  \]
for all $v,w \in \C[D]$. Let $v = \sum_{\gamma \in D} v_{\gamma} e^{\gamma} \in \C[D]^{\Gamma}$. Then $(v,\inv_D(e^{\gamma}))=v_{\gamma}$. This implies $\inv_D(e^{\gamma}) = 0$ if $\gamma \notin I$. Furthermore $\inv_D$ commutes with
%$\rho_D$.
$\rho_D(M)$ for all $M\in  \Gamma$. 

%\medskip

For an isotropic element $\gamma \in D$ we can calculate $\inv_D(e^{\gamma})$ as follows:

\begin{thm} \label{invexplicit}
Let $D$ be a discriminant form of even signature and level dividing $N$ and $\gamma \in I$. Then 
\[ \inv_D(e^{\gamma}) = \sum_{s \in \Gamma(N)\backslash P} \inv_D(e^{\gamma})_s \]
with 
\begin{multline*}
\inv_D(e^{\gamma})_s \, = 
\, \xi(M^{-1}) \, \frac{N}{|\Gamma(N) \backslash \Gamma|} \, \frac{\sqrt{|D_c|}}{\sqrt{|D|}} \, \\
\sum_{\mu \in (a \gamma + D^{c*}) \cap I}  e(d\q_c(\mu -a \gamma)) e(b(\mu,\gamma)) e^{\mu} 
\end{multline*}
if $N=2$ and 
\begin{multline*}
\inv_D(e^{\gamma})_s \, = 
\, \xi(M^{-1}) \, \frac{N}{|\Gamma(N) \backslash \Gamma|} \, \frac{\sqrt{|D_c|}}{\sqrt{|D|}} \, \\
\sum_{\mu \in (a\gamma + D^{c*})\cap I}  e(d\q_c(\mu -a \gamma)) e(b(\mu,\gamma))
\big\{ e^{\mu} + e(\sign(D)/4) e^{-\mu} \big\}
\end{multline*}
if $N \geq 3$ where in both cases $M = \left( \begin{smallmatrix} a & b \\ c & d \end{smallmatrix} \right)$ is any matrix in $\Gamma$ such that $M \infty = s$. 
\end{thm}
{\em Proof:} We can write
\[ \inv_D(e^{\gamma}) = \sum_{s \in \Gamma(N)\backslash P} \inv_D(e^{\gamma})_s \]
with
\[ \inv_D(e^{\gamma})_s = 
\frac{1}{|\Gamma(N) \backslash \Gamma|} \, \sum_{\substack{ M \in \Gamma(N)\backslash \Gamma \\ M \infty = s } } \, \rho_D(M^{-1})e^{\gamma} \, . 
\]
Suppose $N \geq 3$. Let $s \in P$ and $M = \left( \begin{smallmatrix} a & b \\ c & d \end{smallmatrix} \right) \in \Gamma$ such that $M \infty =s$. Then
\begin{align*}
\inv_D(e^{\gamma})_s
&= 
\, \frac{1}{|\Gamma(N) \backslash \Gamma|} \, \sum_{n \in \Z/N\Z} 
\big\{ \rho_D((MT^n)^{-1}) e^{\gamma} +  \rho_D((-MT^n)^{-1}) e^{\gamma} \big\}  \\
&=
\, \frac{1}{|\Gamma(N) \backslash \Gamma|} \, \sum_{n \in \Z/N\Z} 
\rho_D(T^{-n})\rho_D(M^{-1}) \{ e^{\gamma} +  e(\sign(D)/4) e^{-\gamma} \}  \\
&= 
\, \xi(M^{-1}) \, \frac{1}{|\Gamma(N) \backslash \Gamma|} \, \frac{\sqrt{|D_c|}}{\sqrt{|D|}} 
\, \sum_{\mu \in a \gamma + D^{c*}} e(d\q_c(\mu -a \gamma)) e(b(\mu,\gamma)) \\
&
\hspace{4cm} 
\sum_{n \in \Z/N\Z} \rho_D(T^{-n}) \{ e^{\mu} +  e(\sign(D)/4) e^{-\mu} \} \\
&= 
\, \xi(M^{-1}) \, \frac{N}{|\Gamma(N) \backslash \Gamma|} \, \frac{\sqrt{|D_c|}}{\sqrt{|D|}} 
\sum_{\mu \in (a \gamma + D^{c*}) \cap I } \, e(d\q_c(\mu -a \gamma)) e(b(\mu,\gamma)) \\
%\sum_{\substack{ \mu \in I \\ \mu \in a \gamma + D^{c*} }} e(d\q_c(\mu -a \gamma)) e(b(\mu,\gamma)) \\
&
\hspace{6.8cm} 
\{ e^{\mu} +  e(\sign(D)/4) e^{-\mu} \} 
\end{align*}
where we used the above formula for the Weil representation, $\gamma \in I$ and
\[   \sum_{n \in \Z/N\Z} \rho_D(T^{-n}) e^{\mu} = 0  \]
if $\mu \notin I$. For $N=2$ we just drop the second sum. \eop

\medskip

The dimension of the subspace of invariants is given by the trace of the linear map $\inv_D$, i.e.
\[  \dim \C[D]^{\Gamma} = \tr \, \inv_D
  = \sum_{\gamma \in I} (\inv_D(e^{\gamma}),e^{\gamma})
  = \sum_{\gamma \in I} \, \sum_{s \in \Gamma(N)\backslash P}
  (\inv_D(e^{\gamma})_s,e^{\gamma})  \, . \]
The previous theorem implies:

\begin{thm} \label{dimformula}
Let $D$ be a discriminant form of even signature and level dividing $N$. Then 
\[ \dim \C[D]^{\Gamma} = \sum_{s \in \Gamma(N)\backslash P} d_s \]
with 
\[
d_s = \, \xi(M^{-1}) \, \frac{N}{|\Gamma(N) \backslash \Gamma|} \, \frac{\sqrt{|D_c|}}{\sqrt{|D|}} \,
\sum_{\substack{ \gamma \in I \\ (1-a)\gamma \in D^{c*} }} \, e( d \q_c( (1-a)\gamma ) )
\]
if $N=2$ and 
\begin{multline*}
\qquad d_s = 
\, \xi(M^{-1}) \, \frac{N}{|\Gamma(N) \backslash \Gamma|} \, \frac{\sqrt{|D_c|}}{\sqrt{|D|}} \, 
\Bigg\{
\, \sum_{\substack{ \gamma \in I \\ (1-a)\gamma \in D^{c*} }} \, e(d \q_c( (1-a)\gamma ) )  \\
+ e(\sign(D)/4)
   \sum_{\substack{ \gamma \in I \\ (1+a)\gamma \in D^{c*} }} \, e(d\q_c( (1+a)\gamma ) ) 
\, \Bigg\}
\end{multline*}
if $N \geq 3$ where in both cases $M = \left( \begin{smallmatrix} a & b \\ c & d \end{smallmatrix} \right)$ is any matrix in $\Gamma$ such that $M \infty = s$. 
\end{thm}

We describe some properties of the invariants of $\rho_D$ and the projection $\inv_D$.

\begin{prp} \label{abba}
Let $D$ be a discriminant form of even signature and level dividing $N$. Let $v = \sum_{\gamma \in D} v_{\gamma} e^{\gamma} \in \C[D]^{\Gamma}$. Then 
\[ v_{\gamma} = \chi_D(a) v_{a\gamma}  \]
for all $a \in (\Z/N\Z)^*$
%$(a,N)=1$
and $\gamma \in D$. If $\chi_D$ is non-trivial and $H$ is a subgroup of $D$, then
\[  \sum_{\gamma \in H} v_{\gamma} =0 \, . \]
\end{prp}
{\em Proof:} Let $M = \left( \begin{smallmatrix} a & b \\ c & d \end{smallmatrix} \right) \in \Gamma_0(N)$. Then 
\[ v = \rho_D(M) v = \sum_{\gamma \in I} v_{\gamma} \rho_D(M) e^{\gamma}
= \chi_D(a) \sum_{\gamma \in I} v_{\gamma} e^{d\gamma} 
= \chi_D(a) \sum_{\gamma \in I} v_{a\gamma} e^{\gamma}  \, .  \]
For the second statement note that $H$ decomposes into orbits under the action of $(\Z/N\Z)^*$ and
\[ \sum_{a \in (\Z/N\Z)^*} v_{a\gamma} = \sum_{a \in (\Z/N\Z)^*} \chi_D(a) v_{\gamma} 
= v_{\gamma} \sum_{a \in (\Z/N\Z)^*} \chi_D(a) = 0 \, . \]
This proves the proposition. \eop

\begin{prp}  \label{tamtam}
Let $D$ be a discriminant form of even signature with non-trivial $\chi_D$. Then
\[  \inv_D(e^0) = 0 \, . \]
\end{prp}
{\em Proof:} We have $(v,\inv_D(e^0)) = v_0 = 0$ for all $v = \sum_{\gamma \in D} v_{\gamma} e^{\gamma} \in \C[D]^{\Gamma}$. Hence $\inv_D(e^0) = 0$. \eop

\begin{prp}  \label{iggypop}
  Let $D$ be a discriminant form of even signature and $\gamma \in I^{\perp}$.
  % Then $v_{\gamma} = v_0$ for all invariants $v = \sum_{\bt \in D} v_{\bt} e^{\bt} \in \C[D]^{\Gamma}$.
  Then $\inv_D(e^{\gamma}) = \inv_D(e^0)$.
\end{prp}
{\em Proof:} 
Let $v = \sum_{\bt \in D} v_{\bt} e^{\bt} \in \C[D]^{\Gamma}$. Then
%the invariance of $v$ under $S$ implies
\begin{align*}
v_{\gamma} 
&= 
\frac{e(\sign(D)/8)}{\sqrt{|D|}}
                  \sum_{\bt\in D} v_{\bt} \, e((\gamma,\bt))
= 
\frac{e(\sign(D)/8)}{\sqrt{|D|}}
                  \sum_{\bt\in I} v_{\bt} \, e((\gamma,\bt)) \\
&= 
\frac{e(\sign(D)/8)}{\sqrt{|D|}}
                  \sum_{\bt\in I} v_{\bt} 
= 
\frac{e(\sign(D)/8)}{\sqrt{|D|}}
                  \sum_{\bt\in D} v_{\bt} 
= v_0  
\end{align*}
where we used the invariance of $v$ under $S$ in the first and in the last step. It follows $(v,\inv_D(e^{\gamma})) = v_{\gamma} = v_0 = (v,\inv_D(e^0))$. \eop

\medskip
\noindent
Let $D$ be a discriminant form of even signature which contains no non-trivial isotropic element, i.e. $I= \{ 0 \}$. Then $I^{\perp} = D$. If in addition $D$ is non-trivial, then the proposition implies $\inv_D(e^0) = 0$ and $\dim( \C[D]^{\Gamma} ) = 0$.
%\[    \dim( \C[D]^{\Gamma} ) = 0  \, . \]
(Choose $\gamma \in D\backslash \{ 0 \}$. Then $\gamma$ is non-isotropic so that $\inv_D(e^0) = \inv_D(e^{\gamma}) = 0$.) 

\begin{prp}  \label{ironmaiden}
  Let $D$ be a discriminant form of even signature with non-trivial $\chi_D$ and $\gamma \in D$ such that $2\gamma \in I^{\perp}$. Then $\inv_D(e^{\gamma}) = 0$.
\end{prp}
{\em Proof:} 
Let $v = \sum_{\bt \in D} v_{\bt} e^{\bt} \in \C[D]^{\Gamma}$. The invariance of $v$ under $S$ and Proposition \ref{abba} give
\begin{align*}
v_{\gamma} 
&= 
\frac{e(\sign(D)/8)}{\sqrt{|D|}}
                  \sum_{\bt\in D} v_{\bt} \, e((\gamma,\bt))
= 
\frac{e(\sign(D)/8)}{\sqrt{|D|}}
                  \sum_{\bt\in I} v_{\bt} \, e((\gamma,\bt)) \\
&= 
\frac{e(\sign(D)/8)}{\sqrt{|D|}} 
\Bigg( \, \sum_{\substack{ \bt \in I \\ (\bt,\gamma) = 0 \bmod 1}} v_{\bt} 
- \sum_{\substack{ \bt \in I \\ (\bt,\gamma) = 1/2 \bmod 1}} v_{\bt} \, \Bigg) \\
&= 
\frac{e(\sign(D)/8)}{\sqrt{|D|}} 
\Bigg( \, 2 \sum_{\substack{ \bt \in I \\ (\bt,\gamma) = 0 \bmod 1}} v_{\bt} 
- \sum_{\bt \in I} v_{\bt} \, \Bigg) \\
&= 
- v_0 + 2 \, \frac{e(\sign(D)/8)}{\sqrt{|D|}} 
\sum_{\bt \in \gamma^{\perp}} v_{\bt}
%= 0 + 0 = 0 
= 0
\end{align*}
because $0$ and $\gamma^{\perp}$ are subgroups of $D$. \eop
%This implies the proposition. \eop

\begin{prp}  \label{cash}
  Let $D$ be a discriminant form of even signature and level dividing $N$. Suppose $(N,5)=1$ and $\chi_D(5)=-1$. Let $\gamma \in D$ such that $4\gamma \in I^{\perp}$. Then $\inv_D(e^{\gamma}) = 0$.
  
\end{prp}
{\em Proof:} For $v = \sum_{\bt \in D} v_{\bt} e^{\bt} \in \C[D]^{\Gamma}$ we have
\[ v_{\gamma} = 
\frac{e(\sign(D)/8)}{\sqrt{|D|}} \sum_{j=0}^3 e(j/4)
\sum_{\substack{ \bt \in I \\ (\bt,\gamma) = j/4 \bmod 1}} v_{\bt} \, . \]
The sets $\{ \bt \in I \, | \, (\bt,\gamma) = j/4 \bmod 1 \}$ are invariant under multiplication by $5$. On the other hand $v_{5\bt} = \chi_D(5) v_{\bt} = -v_{\bt}$ for all $\bt \in D$ by Proposition \ref{abba}. It follows 
\[  2 \sum_{\substack{ \bt \in I \\ (\bt,\gamma) = j/4 \bmod 1}} v_{\bt} =
    \sum_{\substack{ \bt \in I \\ (\bt,\gamma) = j/4 \bmod 1}} ( v_{\bt} + v_{5\bt} ) = 0 \, . \]
This implies the statement. \eop

\medskip
\noindent
Note that the condition of the proposition is satisfied for example for $2$-adic discriminant forms $D$ such that $|D|$ is not a square.

\section{Discriminant forms of prime level} \label{wims}

In this section we calculate the projection on the subspace of invariants and the dimension of this space explicitly for discriminant forms of prime level.
%in the case that the discriminant form is an elementary abelian $p$-group.

\medskip

We start with the case that $p$ is odd.

\begin{thm} \label{invelementaryodd}
Let $p$ be an odd prime and $D$ a discriminant form of type $p^{\epsilon n}$. Let $\gamma \in I$. Then
\begin{multline*}
\inv_D(e^{\gamma}) = \epsilon \left( \frac{-1}{p} \right)^{n/2}
\frac{1}{p^2-1} \, \frac{1}{p^{(n-2)/2}} \,  
\Bigg\{ \sum_{\mu \in (\gamma^{\perp}\cap I)} p \, e^{\mu} - \sum_{\mu \in I} e^{\mu} \Bigg\}  \\
+ \frac{1}{p^2-1} \, \sum_{a \in (\Z/p\Z)^*} e^{a\gamma} 
\end{multline*}
if $n$ is even and
\begin{multline*}
\inv_D(e^{\gamma}) = \epsilon \left( \frac{-1}{p} \right)^{(n+1)/2} \left( \frac{2}{p} \right) 
\frac{1}{p^2-1} \, \frac{1}{p^{(n-3)/2}} \, \sum_{\mu \in I} \left( \frac{p(\mu,\gamma)}{p} \right) e^{\mu} \\
+ \frac{1}{p^2-1} \, \sum_{a \in (\Z/p\Z)^*} \left( \frac{a}{p} \right) e^{a\gamma} 
\end{multline*}
if $n$ is odd.
\end{thm}
{\em Proof:} The cusps of $\Gamma(p)$ are represented by the pairs $\pm(a,c) \in (\Z/p\Z)^2 \backslash \{ (0,0) \}$. 
If $(c,p)=1$, we can choose any $d \in \Z/p\Z$ and define $b = c^{-1}(ad-1)$ to obtain a matrix $\left( \begin{smallmatrix} a & b \\ c & d \end{smallmatrix} \right) \in \SL_2(\Z/p\Z)$ ($c^{-1}$ denotes the inverse of $c$ modulo $p$). Let $M_s$ be any lift of $\left( \begin{smallmatrix} a & b \\ c & d \end{smallmatrix} \right)$ to $\Gamma$ (recall that the projection $\SL_2(\Z) \to \SL_2(\Z/p\Z)$ is surjective). Then 
\begin{multline*}
  \inv_D(e^{\gamma})_s \, =
  \, \xi(M_s^{-1}) \, \frac{1}{p^2-1} \, \frac{1}{p^{n/2}} \\
\sum_{\mu \in (a\gamma + D^{c*})\cap I}  e(d\q_c(\mu -a \gamma)) e(b(\mu,\gamma))
\big\{ e^{\mu} + e(\sign(D)/4) e^{-\mu} \big\} \, . 
\end{multline*}
Taking $d=0 \! \mod p$ and using the explicit formula for $\xi(M_s^{-1})$ given in \cite{S2} we obtain
\begin{multline*}
\inv_D(e^{\gamma})_s = \, e(\sign(D)/8) \, \left( \frac{-c}{|D|} \right) \, \frac{1}{p^2-1} \, \frac{1}{p^{n/2}} \, \\
\sum_{\mu \in I} \, e(-c^{-1}(\mu,\gamma)) \{ e^{\mu} + e(\sign(D)/4) e^{-\mu} \} \, . 
\end{multline*}
If $c = 0 \! \mod p$, we choose a matrix $\left( \begin{smallmatrix} a & b \\ c & d \end{smallmatrix} \right) \in \SL_2(\Z/p\Z)$ and lift it to a matrix $M_s$ in $\Gamma$. Then 
\[
\inv_D(e^{\gamma})_s = \, \left( \frac{a}{|D|} \right) \, \frac{1}{p^2-1} \, 
\{ e^{a\gamma} + e(\sign(D)/4) e^{-a\gamma} \} \, .
\]
Summing over all cups of $\Gamma(N)$ we get
\begin{multline*}
\inv_D(e^{\gamma}) = \frac{1}{2} \, e(\sign(D)/8) \, \frac{1}{p^2-1} \, \frac{1}{p^{(n-2)/2}} \, 
\sum_{\mu \in I} \{ e^{\mu} + e(\sign(D)/4) e^{-\mu} \}  \\
\sum_{c \in (\Z/p\Z)^*} \left( \frac{c}{|D|} \right) \, e(c(\mu,\gamma)) \\
+ \frac{1}{2} \, \frac{1}{p^2-1} \, 
\sum_{a \in (\Z/p\Z)^*} \left( \frac{a}{|D|} \right) \{ e^{a\gamma} + e(\sign(D)/4) e^{-a\gamma} \}  \, . 
\end{multline*}
If $n$ is even, then $e(\sign(D)/8) = \epsilon \big( \frac{-1}{p} \big)^{n/2}$ (see the proof of Theorem 7.1 in \cite{S1}) and 
\[  
\sum_{c \in (\Z/p\Z)^*} \left( \frac{c}{|D|} \right) \, e(c(\mu,\gamma)) = 
\begin{cases}
\, p-1 &  \text{if $(\mu,\gamma) = 0 \! \! \mod 1$,} \\
\, -1  &  \text{otherwise}
\end{cases}
\]
so that 
\begin{multline*}
\inv_D(e^{\gamma}) = \epsilon \left( \frac{-1}{p} \right)^{n/2}
\frac{1}{p^2-1} \, \frac{1}{p^{(n-2)/2}} \,  
\Bigg\{ \sum_{\mu \in (\gamma^{\perp}\cap I)} p \, e^{\mu} - \sum_{\mu \in I} e^{\mu} \Bigg\}  \\
+ \frac{1}{p^2-1} \, \sum_{a \in (\Z/p\Z)^*} e^{a\gamma} \, .
\end{multline*}
If $n$ is odd, the statement follows from
\[
e(\sign(D)/8) = \epsilon \left( \frac{2}{p} \right)
\begin{cases}
  \, 1                  & \text{if } \, p = 1 \, \bmod 4 \, , \\
  \, (-1)^{(n+1)/2} (-i) & \text{if } \, p = 3 \, \bmod 4  
\end{cases}
\]
and
\[
\sum_{c \in (\Z/p\Z)^*} \left( \frac{c}{p} \right) \, e(c(\mu,\gamma)) 
= \left( \frac{p(\mu,\gamma)}{p} \right) \sqrt{p} \, 
\begin{cases}
\, 1   & \text{if } \, p = 1 \, \bmod 4  \,  , \\
\, i   & \text{if } \, p = 3 \, \bmod 4  
\end{cases}
\]
(see Theorem 1.2.4 in \cite{BEW}). \eop

\medskip
\noindent
Note that for $n=1$ or $n=2$ and $\epsilon \big( \frac{-1}{p} \big) = -1$ we have $I = \{ 0 \}$ and $\inv_D(e^{\gamma}) = 0$ for all $\gamma \in D$. 
The first formula in the theorem extends to discriminant forms of level $2$. 

\begin{thm}
Let $D$ be a discriminant form of type $2_{I\!I}^{\epsilon n}$ with $n$ even and $\gamma \in I$. Then
\[
\inv_D(e^{\gamma}) = \epsilon \, \frac{1}{3} \, \frac{1}{2^{(n-2)/2}} \,  
\Bigg\{ \sum_{\mu \in (\gamma^{\perp}\cap I)} 2 \, e^{\mu} - \sum_{\mu \in I} e^{\mu} \Bigg\}  
+ \frac{1}{3} \, e^{\gamma} \, .
\]
\end{thm}

\noindent
Next we calculate the dimension of the subspace of invariants.

\begin{thm}  \label{dimformulaodd}
Let $p$ be an odd prime and $D$ a discriminant form of type $p^{\epsilon n}$. Then
\[ \dim \C[D]^{\Gamma} = \frac{p^{n-1}-p}{p^2-1} + \epsilon \left( \frac{-1}{p} \right)^{n/2} p^{(n-2)/2} + 1   \]
if $n$ is even and
\[ \dim \C[D]^{\Gamma} = \frac{p^{n-1}-1}{p^2-1} \]
if $n$ is odd.
\end{thm}
{\em Proof:} This can be proved directly by using Theorem \ref{dimformula} or by means of Theorem \ref{invelementaryodd}. We describe the second approach for $n$ even. For $\gamma \in I$ we have
\begin{multline*}
(\inv_D(e^{\gamma}),e^{\gamma}) = 
\epsilon \left( \frac{-1}{p} \right)^{n/2} \frac{1}{p^2-1} \, \frac{1}{p^{(n-2)/2}} \, (p-1) \\
+ \frac{1}{p^2-1} \begin{cases}
\, 1   & \text{if $\gamma \neq 0$}, \\
\, p-1 & \text{if $\gamma = 0$} 
\end{cases}
\end{multline*}
so that
\begin{align*}
\dim \C[D]^{\Gamma} 
&= 
\sum_{\gamma \in I} \, (\inv_D(e^{\gamma}),e^{\gamma}) \\
&=
|I| \, \big\{ \epsilon \left( \frac{-1}{p} \right)^{n/2} \frac{1}{p^2-1} \, \frac{1}{p^{(n-2)/2}} \, (p-1) \big\} \\
& 
\hspace*{2.8cm} + \, \frac{1}{p^2-1} \, \big\{ |I\backslash \{ 0\}| + (p-1) \big\} \\
&= 
\frac{p^{n-1}-p}{p^2-1} + \epsilon \left( \frac{-1}{p} \right)^{n/2} p^{(n-2)/2} + 1 
\end{align*}
by Proposition \ref{numberofelementspodd}. \eop

\medskip
\noindent
We describe an example. If $D$ is of type $p^{\epsilon 2}$ with $\epsilon = \big( \frac{-1}{p} \big)$, the subspace of invariants has dimension $\dim \C[D]^{\Gamma} = 2$. The discriminant form $D$ is generated by two isotropic elements $\gamma_1, \gamma_2$ such that $(\gamma_1,\gamma_2) = 1/p \! \mod 1$. We have
\[   \inv_D(e^0) = \frac{1}{p+1} \, 
\Bigg\{ e^0 + \sum_{\mu \in I} \, e^{\mu} \Bigg\}  \]
and
\[   \inv_D(e^{\gamma_i}) = 
\frac{1}{p-1} \, \sum_{\mu \in \la \gamma_i \ra} e^{\mu} 
- \frac{1}{p^2-1} \Bigg\{ e^0 + \sum_{\mu \in I} \, e^{\mu} \Bigg\}  \, . \]
This implies that $\C[D]^{\Gamma}$ is generated by the elements $\sum_{\mu \in \la \gamma_i \ra} e^{\mu}$, $i = 1,2$,
% and $\sum_{\mu \in \la \gamma_2 \ra} e^{\mu}$,
which is a special case of Skoruppa's result.

\medskip

As for odd primes we can prove

\begin{thm} \label{dimformulaeven}
Let $D$ be a discriminant form of type $2_{I\!I}^{\epsilon n}$ with $n$ even. Then
\[   \dim \C[D]^{\Gamma} = \frac{2^{n-1}+1}{3} + \epsilon \, 2^{(n-2)/2} \, .   \]
\end{thm}

The dimension formulas in Theorems \ref{dimformulaodd} and \ref{dimformulaeven}
have also been found by Zemel using a slightly different approach (see Theorem 5.6 in \cite{Z}). We also remark that the numerical values of $\dim (\C[D]^{\Gamma})$ for some of the above discriminant forms and others have been determined by Skoruppa and Ehlen (see Section 6 in \cite{ES}).

\begin{cor}  \label{invp3}
  Let $p$ be an odd prime and $D$ a discriminant form of type $p^{\epsilon 3}$ with $\epsilon = \pm 1$. Choose $\gamma \in I\backslash \{ 0\}$. For $j \in \Z/p\Z$ define 
\[  
M(\gamma)_j = \{ \, \mu \in I\backslash \{ 0 \} \, | \, (\mu,\gamma) = j/p \! \! \mod 1 \, \}  \, . 
\]
Let
\[  M(\gamma)^+ = 
\bigcup_{\substack{ j \in (\Z/p\Z)^* \\ \varepsilon \chi_D(j)=+1} } M(\gamma)_j 
\cup \bigcup_{\substack{ j \in (\Z/p\Z)^* \\ \chi_D(j)=+1} } \{ j\gamma \} \]
where $\varepsilon = \epsilon \big( \frac{2}{p} \big)$ and analogously $M(\gamma)^-$.
  %$I^+ = \{ \, \mu \in I \, | \, \epsilon \big( \frac{2}{p} \big) \big( \frac{p(\mu,\gamma)}{p} \big) = +1 \, \} \cup \{ \, a \gamma \, | \, a \in (\Z/p\Z)^* , \big( \frac{a}{p} \big) = +1 \}$ 
%\[   M(\gamma)^+ =  \Big\{ \, \mu \in I \, \Big| \, \epsilon \left( \frac{2}{p} \right) \left( \frac{p(\mu,\gamma)}{p} \right) = +1 \, \Big\}  \cup  \Big\{ \, a \gamma \, \Big| \left( \frac{a}{p} \right) = +1 \, \Big\}  \]
  Then $\C[D]^{\Gamma}$ is $1$-dimensional and spanned by 
\[  \sum_{\mu \in M(\gamma)^+} e^{\mu} - \sum_{\mu \in M(\gamma)^-} e^{\mu}  \, . \]

If $D$ is of type $p^{-4}$ where $p$ is an odd prime, then $\C[D]^{\Gamma}$ is $1$-dimensional and spanned by 
\[  p e^0 - \sum_{\mu \in I} e^{\mu}  \, . \]
The same result holds for $D$ of type $2_{I\!I}^{-4}$.
\end{cor}
{\em Proof:} In the first case $\C[D]^{\Gamma}$ is spanned by $\inv_D(e^{\gamma})$ for any $\gamma \in I\backslash \{ 0\}$ and in the second case by $\inv_D(e^0)$. \eop 

\medskip
\noindent
The decomposition $I\backslash \{ 0\} = M(\gamma)^+ \cup M(\gamma)^-$ is independent of the choice of $\gamma \in I\backslash \{ 0\}$ and is equal to the decomposition of $I\backslash \{ 0\}$ under the action of the spinor kernel of $\text{SO}(D)$. The size of $M(\gamma)^{\pm}$ is $(p^2-1)/2$.

\section{Some $2$-adic exercises} \label{zwei}

We study some $2$-adic discriminant forms which will play an important role in our main theorem.

\medskip

Let $D$ be a discriminant form of type $2_t^{\epsilon n}$. Then $D^{2*}$ contains a single element which we denote by $x_2$. The signature of $D$ is even if and only if $n$ is even. In this case the matrix $Z = -1 \in \Gamma$ acts as $\rho_D(Z) e^{\gamma} = e(t/4) e^{\gamma}$ for all $\gamma \in D$ so that there are no non-trivial invariants if $t = 2 \! \mod 4$.

\begin{prp}  \label{pba}
Let $D$ be a discriminant form of type $2_t^{\epsilon n}$ with $n$ even and $t = 0 \! \mod 4$. Then 
\[
\inv_D(e^{\gamma}) = 
\frac{1}{6} e^{\gamma} + \frac{1}{6} e^{\gamma + x_2} 
+ \epsilon \, (-1)^{t/4} \, \frac{1}{6} \, \frac{1}{2^{(n-4)/2}} \,  
\Bigg\{ \sum_{\mu \in (\gamma^{\perp}\cap I)} 2 \, e^{\mu} - \sum_{\mu \in I} e^{\mu} \Bigg\}  
\]
for $\gamma \in I$ and 
\[   \dim \C[D]^{\Gamma} = \frac{2^{n-3} + 1}{3} + \epsilon(-1)^{t/4}2^{(n-4)/2}  \, . \]
\end{prp}
The proof is similar to the proof of the next theorem. We therefore omit it.

\medskip
\noindent
We describe two examples. If $D$ is of type $2_0^{+2}$, then $\C[D]^{\Gamma}$ is $1$-dimensional and spanned by $\inv_D(e^{0}) = \inv_D(e^{x_2}) = (e^{0} + e^{x_2})/2$. If $D$ is of type
 $2_0^{-4} \cong 2_4^{+4}$,
%$2_0^{-4} = 2_4^{+4}$,
then $\C[D]^{\Gamma}$ is trivial.

\medskip

Let $D$ be a discriminant form of type $2_t^{\epsilon n} 4_{I\!I}^{+2}$. Then $D$ has even signature if and only if $n$ is even. 

\begin{prp}
Let $D$ be a discriminant form of type $2_t^{\epsilon n} 4_{I\!I}^{+2}$ with $n$ even. Then
\begin{multline*}
\inv_D(e^{\gamma}) = 
\frac{1}{12} \{ e^{\gamma} + e(t/4) e^{-\gamma} \} \\
+ \frac{1}{24} \sum_{\mu \in (\gamma + D^{2*}) \cap I} 
                e(\q_2(\mu - \gamma)) \{ e^{\mu} + e(t/4) e^{-\mu} \} \\
+ \epsilon \, e(3t/8) \, \frac{1}{12} \, \frac{1}{2^{n/2}} \sum_{\mu \in I} 
                e(-(\mu,\gamma)) \{ e^{\mu} + e(t/4) e^{-\mu} \} 
\end{multline*}
for $\gamma \in I$ and 
\[   \dim \C[D]^{\Gamma} =  
\frac{1}{12} \, |I| \big\{ 1 + \epsilon \, e(3t/8) \frac{1}{2^{n/2}} (1 + e(t/4)) \big\} 
   + \frac{1}{12}\, e(t/4) \, |I_2| \]
 with
 \begin{align*}
   |I|   &= 2^{n+2} + \epsilon \, 2^{(n+2)/2} \, \delta(t/4) \left( \frac{t-1}{2} \right)    \\
   |I_2| &= 2^n    + \epsilon \, 2^{(n+2)/2} \, \delta(t/4) \left( \frac{t-1}{2} \right) .
 \end{align*}
\end{prp}
{\em Proof:} First note that $e(\sign(D)/8) = \gamma_2(D) = \epsilon \, e(t/8)$. The group $\Gamma(4)$ has $6$ cusps $s$, which can be represented by $1/4, 1/2$ and $a/1$ with $a=0,1,2,3$. Choosing matrices $M_s$ as
$\left( \begin{smallmatrix} 1 & 0  \\ 4 & 1 \end{smallmatrix} \right),
 \left( \begin{smallmatrix} 1 & 0  \\ 2 & 1 \end{smallmatrix} \right)$ and 
$\left( \begin{smallmatrix} a & -1 \\ 1 & 0 \end{smallmatrix} \right)$
we find 
\[  \inv_D(e^{\gamma})_{1/4} = \frac{1}{12} \{ e^{\gamma} + e(t/4) e^{-\gamma} \}  \]
and 
\[  \inv_D(e^{\gamma})_{1/2} = \frac{1}{24} \sum_{\mu \in (\gamma + D^{2*}) \cap I} 
                e(\q_2(\mu - \gamma)) \{ e^{\mu} + e(t/4) e^{-\mu} \}   \, .   \]
We remark that in the last sum $e(\q_2(\mu - \gamma)) = \pm 1$. For $M_s = \left( \begin{smallmatrix} a & -1 \\ 1 & 0 \end{smallmatrix} \right)$ we have $\xi(M_s^{-1}) = \epsilon \, e(3t/8)$ and $a\gamma + D^{1*} = D$ so that
\[  \inv_D(e^{\gamma})_{a/1} = 
\epsilon \, e(3t/8) \, \frac{1}{48} \, \frac{1}{2^{n/2}} \sum_{\mu \in I} 
                e(-(\mu,\gamma)) \{ e^{\mu} + e(t/4) e^{-\mu} \} \, .  \]
This implies the formula for $\inv_D(e^{\gamma})$.

Next we calculate the dimension of the fixed point subspace. For $\gamma \in I$ we have
\begin{align*}
(\inv_D(e^{\gamma})_{1/4},e^{\gamma})
  &= \frac{1}{12} + \frac{1}{12} \, e(t/4) \,
\begin{cases}
\, 1   & \text{if $2\gamma = 0$,} \\
% \, 0   & \text{if $2\gamma \neq 0$ ,}
\, 0   & \text{otherwise}
\end{cases}    \\
(\inv_D(e^{\gamma})_{1/2},e^{\gamma})
& = 0  \\
(\inv_D(e^{\gamma})_{a/1},e^{\gamma}) 
&= \epsilon \, e(3t/8) \, \frac{1}{48} \, \frac{1}{2^{n/2}} (1 + e(t/4))
\end{align*}
so that
\begin{align*}
 \dim \C[D]^{\Gamma}
&= \sum_{\gamma \in I} (\inv_D(e^{\gamma}),e^{\gamma})  \\
&= \frac{1}{12} \, |I| \big\{ 1 + \epsilon \, e(3t/8) \, \frac{1}{2^{n/2}} (1 + e(t/4)) \big\} 
   + \frac{1}{12}\, e(t/4) \, |I_2|
\end{align*} 
where $I_2 = I \cap D_2$. The cardinalities of $I$ and $I_2$ can be determined with Proposition \ref{numberofelements2odd}. \eop

\medskip
\noindent
Note that if $t = 2 \! \mod 4$ and $2\gamma = 0$ then $\inv_D(e^{\gamma}) = 0$. This also follows from the formula for the action of
%$Z=-1$.
$Z = \left( \begin{smallmatrix} -1 & 0  \\ 0 & -1 \end{smallmatrix} \right)$.

\medskip

Now we consider the case that $D$ is of type $2_t^{+2} 4_{I\!I}^{+2}$ with $t = 2 \! \mod 4$. Then $\sign(D) = t \! \mod 8$. The set $I\backslash I_2$ has cardinality $16 - 4 = 12$ and $\Orth(D)$ acts transitively on it. Let $\gamma \in I\backslash I_2$. For $j \in \Z/4\Z$ we define 
\[  
M(\gamma)_j = \{ \, \mu \in I\backslash I_2 \, | \, (\mu,\gamma) = j/4 \! \! \mod 1 \, \}  \, . 
\]
Then
\[  
| M(\gamma)_j | =    
\begin{cases}
\, 4  & \text{if $j$ is odd,} \\
\, 2  & \text{if $j$ is even.} 
\end{cases}
\]
We have $M(\gamma)_0 = \{ \pm \gamma \}$. There is a unique element $\mu \in D^{2*}$ such that $\q_2(\mu) = 0 \! \mod 1$ and $(\mu,\gamma) = 1/2 \! \mod 1$. Define $\al = \mu + \gamma$. Then $M(\gamma)_2 = \{ \pm \al \}$. Let
\[  M(\gamma)^+ =  M(\gamma)_j \, \cup \, \{ +\al \} \, \cup \, \{ +\gamma \}  \]
with $j \in (\Z/4\Z)^*$ such that $\varepsilon \chi_d(j) = +1$ and
\[  M(\gamma)^- =  M(\gamma)_j \, \cup \, \{ -\al \} \, \cup \, \{ -\gamma \}  \]
with $j \in (\Z/4\Z)^*$ such that $\varepsilon \chi_D(j) = -1$ where in both cases
\[ \varepsilon =
\begin{cases}
\, 1   & \text{if } \, t = 6 \, \bmod 8 \, ,  \\
\, -1  & \text{if } \, t = 2 \, \bmod 8 \, .
\end{cases}
\]
%and analogously $M(\gamma)^-$.

\begin{prp} \label{inv24}
Let $D$ be a discriminant form of type $2_t^{+2} 4_{I\!I}^{+2}$ with $t = 2 \! \mod 4$. 
Then the subspace of invariants $\C[D]^{\Gamma}$ is $1$-dimensional and spanned by
\[  \sum_{\mu \in M(\gamma)^+} e^{\mu} - \sum_{\mu \in M(\gamma)^-} e^{\mu}  \]
where $\gamma$ is any element in $I \backslash I_2$.
\end{prp}
{\em Proof:}
By the previous proposition 
\[  \dim \C[D]^{\Gamma} = \frac{1}{12} ( |I| - |I_2| ) =  \frac{1}{12} ( 16 - 4 ) =  1 \, . \]
For $\gamma \in I \backslash I_2$ we have
\begin{multline*}
\inv_D(e^{\gamma}) = 
\frac{1}{12} \{ e^{\gamma} - e^{-\gamma} \} \\
+ \frac{1}{24} \sum_{\mu \in (\gamma + D^{2*}) \cap I} 
                e(\q_2(\mu - \gamma)) \{ e^{\mu} - e^{-\mu} \} \\
+ e(3t/8) \, \frac{1}{24} 
  \sum_{\substack{ \mu \in I \backslash I_2 \\ (\mu,\gamma) = \pm 1/4 } } e(-(\mu,\gamma)) \{ e^{\mu} - e^{-\mu} \} \, .
\end{multline*}
The sum is supported on $I \backslash I_2$ (see Proposition \ref{ironmaiden}). We easily check that 
\[ \inv_D(e^{\gamma}) = 
\frac{1}{12} \, \Bigg\{ \sum_{\mu \in M(\gamma)^+} e^{\mu} - \sum_{\mu \in M(\gamma)^-} e^{\mu}  \Bigg\} \, . \]
This proves the proposition. \eop

\medskip
\noindent
Note that the decomposition $I \backslash I_2 = M(\gamma)^+ \cup M(\gamma)^-$ is independent of the choice of $\gamma \in I \backslash I_2$.
The size of $M(\gamma)^{\pm}$ is $4+1+1=6=12/2$. 

\medskip

We remark that every discriminant form $D$ of level $4$, exponent $4$, order $4^3$ and signature $t = 2 \! \mod 4$ is isomorphic to $2_t^{+2} 4_{I\!I}^{+2}$.

\medskip

Next we consider a discriminant form $D$ of type $2_1^{+1} 4_t^{\epsilon} \, 8_{I\!I}^{+2}$ with $t = 1 \! \mod 2$ and $\epsilon = \left( \frac{t}{2} \right)$. Then $\sign(D) = 1 + t \! \mod 8$. Recall that $I_4 = I \cap D_4$.

\begin{prp} \label{ref}
We have $|I|=64$ and $|I_4|=16$.
\end{prp}
{\em Proof:}
The partition function of $8_{I\!I}^{+2}$ is given by
\[  f_{8_{I\!I}^{+2}}(x) = \sum_{\gamma \in 8_{I\!I}^{+2}} x^{8 \! \q(\gamma)} =
20 + 4(x+x^3+x^5+x^7) + 8(x^2+x^6) + 12x^4 \]
where we have chosen $\q(\gamma) \in [0,1)$. Multiplying this polynomial with the polynomials $f_{2_1^{+1}}(x) = 1 + x^2$ and $f_{4_t^{\epsilon}}(x) = 1 + 2 x^t + x^4$ we can easily derive the first statement. The second follows analogously. \eop

\begin{prp}  \label{oneorbit}
The group $\Orth(D)$ acts transitively on $I\backslash I_4$.
\end{prp}
{\em Proof:} 
Let $\gamma \in I\backslash I_4$. Then there is an element $\bt \in D$ such that $(\gamma,\bt) = 1/8 \! \mod 1$. Define $\mu = \bt - a \gamma$ where $a = 8 \q(\bt) \! \mod 8$. Then $\la \gamma, \mu \ra$ is a discriminant form of type $8_{I\!I}^{+2}$. The orthogonal complement $\la \gamma, \mu \ra^{\perp}$ is a discriminant form of type $2_{t_2}^{\epsilon_2} 4_{t_4}^{\epsilon_4}$. Up to isomorphism there are exactly $4$ such forms namely the forms of type $2_1^{+1} 4_{t_4}^{\epsilon_4}$ with $t_4$ odd and $\epsilon_4 = \left( \frac{t_4}{2} \right)$. The signature of such a form is $1 + t_4 \! \mod 8$. Hence each element $\gamma$ in $I\backslash I_4$ gives rise to a Jordan decomposition $2_1^{+1} 4_t^{\epsilon} \, 8_{I\!I}^{+2}$. This implies that all elements in $I\backslash I_4$ are conjugate under $\Orth(D)$. \eop 

\begin{prp} \label{vi4}
% Let $D$ be a discriminant form as in Proposition \ref{ref}.
% Let $D$ be a discriminant form of type $2_1^{+1} 4_t^{\epsilon} \, 8_{I\!I}^{+2}$ with $t = 1 \! \mod 2$ and $\epsilon = \left( \frac{t}{2} \right)$.
  Let $\gamma \in I_4$. Then $\inv_D(e^{\gamma}) = 0$.
\end{prp}
{\em Proof:}
We have $4 \gamma = 0 \in I^{\perp}$ so that $\inv_D(e^{\gamma})=0$ by Proposition \ref{cash}. \eop

\begin{prp} \label{firstformula}
For $\gamma \in I$ we have
\begin{align*}
\MoveEqLeft{\inv_D(e^{\gamma}) = }  \\
& e(-\sign(D)/8) \, \frac{1}{48} \, \frac{1}{2\sqrt{2}} 
    \sum_{\mu \in I} e(-(\mu,\gamma)) (1-e(-4(\mu,\gamma))) \\
&  \hspace*{7cm} \{ e^{\mu} + e(\sign(D)/4) e^{-\mu} \} \\
& + \epsilon \, e(-t/8) \, \frac{1}{48} \, \frac{1}{4 \sqrt{2}} 
\sum_{\substack{ a \in \Z/8\Z \\ a = 1 \bmod 2} } \,
\sum_{\mu \in (a\gamma + D^{2*}) \cap I} 
e(\q_2(\mu - a \gamma)) \, e(\tfrac{a-1}{2}(\mu,\gamma))  \\
& \hspace*{7cm}  \{ e^{\mu} + e(\sign(D)/4) e^{-\mu} \}  \\
& + \frac{1}{48} \, \frac{1}{2} \, 
\sum_{\substack{ a \in \Z/8\Z \\ a = 1 \bmod 4} } \,
\sum_{\mu \in (a\gamma + D^{4*}) \cap I} 
e(\q_4(\mu - a \gamma)) \, e(\tfrac{(a-1)}{4}(\mu,\gamma))  \\
&  \hspace*{7cm}  \{ e^{\mu} + e(\sign(D)/4) e^{-\mu} \} \\
& + \frac{1}{48} \, \sum_{\substack{ a \in \Z/8\Z \\ a = 1 \bmod 2} } \chi_D(a) \, e^{a \gamma} .
\end{align*}
\end{prp}
{\em Proof:}
The group $\Gamma(8)$ has $24$ cusps. 
There are $16$ cusps $s=a/c \in \Q$, $(a,c)=1$ with $(c,8)=1$. For such a cusp we can choose a matrix $M_s = \left( \begin{smallmatrix} a & b \\ c & d \end{smallmatrix} \right) \in \Gamma$ with $d = 0 \! \mod 8$. Then  $b = -c \! \mod 8$ and $\xi(M_s^{-1}) = \left( \frac{c}{2} \right) e(-c\sign(D)/8)$ so that 
\begin{multline*}
\inv_D(e^{\gamma})_s = 
\left( \frac{c}{2} \right) \, e(-c\sign(D)/8) \, \frac{1}{48} \, \frac{1}{16 \sqrt{2}}  \\
\sum_{\mu \in I} e(-c(\mu,\gamma)) \{ e^{\mu} + e(\sign(D)/4) e^{-\mu} \} \, .
\end{multline*}
There are $4$ cusps $s=a/c \in \Q$, $(a,c)=1$ with $(c,8)=2$. We can choose a matrix $M_s = \left( \begin{smallmatrix} a & b \\ c & d \end{smallmatrix} \right) \in \Gamma$ such that $d = 1 \! \mod 16$. Then 
%$b = \frac{c}{2} \frac{a-1}{2} \! \mod 8$  
$b = c(a-1)/4 \! \mod 8$ and $\xi(M_s^{-1}) = \epsilon \, e(-t/8)$. It follows
\begin{multline*}
\inv_D(e^{\gamma})_s = 
\epsilon \, e(-t/8) \, \frac{1}{48} \, \frac{1}{4 \sqrt{2}}  \\
\sum_{\mu \in (a\gamma + D^{2*}) \cap I} e(\tfrac{c}{2} \q_2(\mu - a \gamma)) \,
e(\tfrac{c(a-1)}{4}(\mu,\gamma)) \, \{ e^{\mu} + e(\sign(D)/4) e^{-\mu} \}  \, .
\end{multline*}
There are $2$ cusps $s=a/c \in \Q$, $(a,c)=1$ with $(c,8)=4$. We choose a representative $s=a/c$ with $a = 1 \! \mod 4$. Then there is a matrix
%$M_s = \left( \begin{smallmatrix} a & b \\ c & d \end{smallmatrix} \right) \in \Gamma$
$M_s = \left( \begin{smallmatrix} a & b \\ c & d \end{smallmatrix} \right)$ in $\Gamma$
such that $d = 1 \! \mod 32$. We find $b = c(a-1)/16 \! \mod 8$ and $\xi(M_s^{-1}) = 1$ so that
\begin{multline*}
\inv_D(e^{\gamma})_s = 
\frac{1}{48} \, \frac{1}{2}  \\
\sum_{\mu \in (a\gamma + D^{4*}) \cap I} e(\tfrac{c}{4} \q_4(\mu - a\gamma)) \,
e(\tfrac{c(a-1)}{16}(\mu,\gamma)) \, \{ e^{\mu} + e(\sign(D)/4) e^{-\mu} \}  \, .
\end{multline*}
Finally there are $2$ cusps $s=a/c \in \Q$, $(a,c)=1$ with $(c,8)=8$. We choose a matrix $M_s = \left( \begin{smallmatrix} a & b \\ c & d \end{smallmatrix} \right) \in \Gamma$. Then $a = d \! \mod 8$ and
the root of unity $\xi(M_s^{-1})$ is given by
$\xi(M_s^{-1}) = \left( \frac{a}{2} \right) e((1-a)\sign(D)/8)$ so that 
\[  \inv_D(e^{\gamma})_s = \left( \frac{a}{2} \right) e((1-a)\sign(D)/8) \, 
\frac{1}{48} \, \{ e^{a\gamma} + e(\sign(D)/4) e^{-a\gamma} \}  \, . \]
Putting the contributions of the cusps together we obtain the given formula. \eop

\begin{prp}
Let $D$ be a discriminant form of type $2_1^{+1} 4_t^{\epsilon} \, 8_{I\!I}^{+2}$ with $t = 1 \! \mod 2$ and $\epsilon = \left( \frac{t}{2} \right)$.
  Then $\C[D]^{\Gamma}$ is $1$-dimensional. 
\end{prp}
{\em Proof:}  
Proposition \ref{vi4} implies  
\begin{multline*} 
\dim \C[D]^{\Gamma}
= \sum_{\gamma \in I} (\inv_D(e^{\gamma}),e^{\gamma})  
= \sum_{\gamma \in I \backslash I_4} (\inv_D(e^{\gamma}),e^{\gamma})   \\ 
= \sum_{\gamma \in I \backslash I_4} \, \sum_{s \in \Gamma(N)\backslash P} 
(\inv_D(e^{\gamma})_s,e^{\gamma})  \, . 
\end{multline*} 
The cusps $s=a/c$ with $(c,8)=1$ do not contribute to the last sum because $1-e(-4(\mu,\gamma))=0$ for $\mu = \pm \gamma$. For $\gamma \in I$ we have $\pm \gamma \notin a \gamma + D^{2*}$ because $\q(x_2) = 1/4 \! \mod 1$ and analogously $\pm \gamma \notin a \gamma + D^{4*}$ because $\q(x_4) = 1/2 \! \mod 1$. Hence the only contribution to the last sum comes from the cusp $1/8$. It follows
\[ \dim \C[D]^{\Gamma} = \frac{1}{48} \, \sum_{\gamma \in I \backslash I_4} 1
= \frac{1}{48} ( |I| - |I_4| ) =  \frac{1}{48} (64 - 16) = 1 \, . \]
This proves the proposition. \eop

\medskip
Finally we show that the generator of $\C[D]^{\Gamma}$ can be written analogously to the cases $p^{\epsilon 3}$ and $2_t^{+2} 4_{I\!I}^{+2}$ (see Corollary \ref{invp3} and Proposition \ref{inv24}). Fix a Jordan decomposition $2_1^{+1} 4_t^{\epsilon} \, 8_{I\!I}^{+2}$ with $t = 1 \! \mod 2$ and $\epsilon = \left( \frac{t}{2} \right)$ of $D$.
Let $\gamma \in I\backslash I_4$.
%Recall that by Proposition \ref{oneorbit} the elements in $I\backslash I_4$ are conjugate under $\Orth(D)$.
For $j \in \Z/8\Z$ we define 
\[  
M(\gamma)_j = \{ \, \mu \in I\backslash I_4 \, | \, (\mu,\gamma) = j/8 \! \! \mod 1 \, \}  \, . 
\]
Then $aM(\gamma)_j = M(\gamma)_{aj}$ for all $a \in (\Z/8\Z)^*$ and
\[  
| M(\gamma)_j | =    
\begin{cases}
\, 8  & \text{if $j$ is odd,} \\
\, 4  & \text{if $j$ is even.} 
\end{cases}
\]
We describe the sets $M(\gamma)_j$ explicitly for even $j$. We have 
\[  M(\gamma)_0 = \{ j \gamma \, | \, j \in (\Z/8\Z)^* \} \, .  \]
There is a unique element $\mu \in D^{4*}$ such that $\q_4(\mu) = 0 \! \mod 1$ and $(\mu,\gamma) = 1/2 \! \mod 1$. Define $\al_4 = \mu + \gamma$. Then
\[  M(\gamma)_4 = \{ j \al_4 \, | \, j \in (\Z/8\Z)^* \} \, .  \]
Finally there are exactly two elements $\mu_i \in D^{2*}$, $i = 1,2$ such that $\q_2(\mu_i) = t/4 \! \mod 1$ and $(\mu_i,\gamma) = 1/4 \! \mod 1$. Define $\al_i = \mu_i + \gamma$. Then
\[  M(\gamma)_2 = \{ \al_1, 5\al_1, \al_2, 5\al_2  \}
\quad \text{and} \quad
M(\gamma)_6 = \{ 3\al_1, 7\al_1, 3\al_2, 7\al_2  \} \, . \]
These statements can be proved by choosing generators of $2_1^{+1} 4_t^{\epsilon} \, 8_{I\!I}^{+2}$ and assuming that $\gamma$ is one of the two isotropic generators of $8_{I\!I}^{+2}$. (Recall that by Proposition \ref{oneorbit} the elements in $I\backslash I_4$ are conjugate under $\Orth(D)$.)
We decompose $I \backslash I_4 = M(\gamma)^+ \cup M(\gamma)^-$ with
\[  M(\gamma)^+ = 
\bigcup_{\substack{ j \in (\Z/8\Z)^* \\ \varepsilon \chi_D(j)=+1} } M(\gamma)_j 
\cup \bigcup_{\substack{ j \in (\Z/8\Z)^* \\ \chi_D(j)=+1} } \{ j\al_1, j\al_2, j\al_4, j\gamma \} \]
and
\[  M(\gamma)^- = 
\bigcup_{\substack{ j \in (\Z/8\Z)^* \\ \varepsilon \chi_D(j)=+1} } M(\gamma)_j 
\cup \bigcup_{\substack{ j \in (\Z/8\Z)^* \\ \chi_D(j)=-1} } \{ j\al_1, j\al_2, j\al_4, j\gamma \} \]
where
\[
  \varepsilon =
\begin{cases}
\, 1   & \text{if } \, t = 5 \text{ or } 7 \, \bmod 8 \, ,  \\
\, -1  & \text{if } \, t = 1 \text{ or } 3 \, \bmod 8 \, .
\end{cases}
\]

\begin{prp} \label{daswarhart}
  Let $D$ be a discriminant form of type $2_1^{+1} 4_t^{\epsilon} \, 8_{I\!I}^{+2}$ with $t = 1 \! \mod 2$ and $\epsilon = \left( \frac{t}{2} \right)$. Then
  %the subspace of invariants
  $\C[D]^{\Gamma}$ is spanned by
\[  \sum_{\mu \in M(\gamma)^+} e^{\mu} - \sum_{\mu \in M(\gamma)^-} e^{\mu}  \]
where $\gamma$ is any element in $I \backslash I_4$.
\end{prp}
{\em Proof:}
Let $\gamma \in I\backslash I_4$. We write $\inv_D(e^{\gamma}) = \sum_{\mu \in I} c_{\mu} e^{\mu}$. Then $c_{\mu} = 0$ for $\mu \in I_4$ by Proposition \ref{vi4}.
%This follows from $\inv_D(e^{\mu})=0$ so that $0=(\inv_D(e^{\gamma}),\inv_D(e^{\mu}))=c_{\mu}$ or directly from the formula in Proposition \ref{firstformula}.
Now we consider the individual sums in Proposition \ref{firstformula}. The first sum extends over 
$\bigcup_{j \in (\Z/8\Z)^*} M(\gamma)_j$, the second over $M(\gamma)_2 \cup M(\gamma)_6$, the third over $M(\gamma)_4$ and the last sum over $M(\gamma)_0$.
For the first sum we find
\begin{align*}
&
\hspace*{-1cm}
\sum_{\mu \in I \backslash I_4} e(-(\mu,\gamma)) (1-e(-4(\mu,\gamma))) 
\{ e^{\mu} + e(\sign(D)/4) e^{-\mu} \}  \\
&= 
2 \, \sum_{j \in (\Z/8\Z)^*} \sum_{\mu \in M(\gamma)_j}
\{ e(-j/8) + e(\sign(D)/4) e(j/8) \} e^{\mu} \\ 
&=
2 \sqrt{2} \,
\sum_{j \in (\Z/8\Z)^*} \chi_D(j) \sum_{\mu \in M(\gamma)_j} e^{\mu} \,
\begin{cases}
\, 1   & \text{if } \,  t = 3 \, \bmod 4, \\
\, -i  & \text{if } \,  t = 1 \, \bmod 4
\end{cases}
\end{align*}         
so that 
\begin{multline*}
e(-\sign(D)/8) \, \frac{1}{48} \, \frac{1}{2\sqrt{2}} 
\sum_{\mu \in I \backslash I_4} e(-(\mu,\gamma)) (1-e(-4(\mu,\gamma))) \\
\hspace*{3.5cm}
\{ e^{\mu} + e(\sign(D)/4) e^{-\mu} \} \\
= \varepsilon \, \frac{1}{48} \, 
\sum_{j \in (\Z/8\Z)^*} \chi_D(j) \sum_{\mu \in M(\gamma)_j} e^{\mu} \, .  
\end{multline*}
We calculate the second sum as  
\begin{align*}
&
\hspace*{-1cm}
\epsilon \, e(-t/8) 
\sum_{\substack{ a \in \Z/8\Z \\ a = 1 \bmod 2} } \,
\sum_{\mu \in (a\gamma + D^{2*}) \cap (I\backslash I_4)} 
e(\q_2(\mu - a \gamma)) \, e(\tfrac{a-1}{2}(\mu,\gamma)) \\
& 
\hspace*{3cm} 
\{ e^{\mu} + e(\sign(D)/4) e^{-\mu} \}  \\
&= 
\epsilon \, e(-t/8) \sum_{\mu\in M(\gamma)_2} \, \sum_{a \in (\Z/8\Z)^*} 
\{ e((a-1)/8) e(\q_2(\mu - a \gamma)) + \\
&
\hspace*{3cm} 
e(\sign(D)/4) e(3(a-1)/8) e(\q_2(-\mu - a \gamma)) \} e^{\mu} \\ 
&
\hspace*{4mm}
+ \epsilon \, e(-t/8) \sum_{\mu\in M(\gamma)_6} \, \sum_{a \in (\Z/8\Z)^*} 
\{ e(3(a-1)/8) e(\q_2(\mu - a \gamma)) + \\
&
\hspace*{3cm} 
e(\sign(D)/4) e((a-1)/8) e(\q_2(-\mu - a \gamma)) \} e^{\mu} \\
&
= 4 \sqrt{2} \, \Bigg\{
\sum_{j \in (\Z/8\Z)^*} \chi_D(j) \, e^{j \al_1} + 
\sum_{j \in (\Z/8\Z)^*} \chi_D(j) \, e^{j \al_2}
\Bigg\} \, .
\end{align*}
Finally we consider the third sum. We easily see that
\begin{multline*}
\sum_{\substack{ a \in \Z/8\Z \\ a = 1 \bmod 4} } \;
\sum_{\mu \in (a\gamma + D^{4*}) \cap (I\backslash I_4)} 
e(\q_4(\mu - \gamma)) \, e(\tfrac{(a-1)}{4}(\mu,\gamma)) \\
\hspace*{4cm}  \{ e^{\mu} + e(\sign(D)/4) e^{-\mu} \}  \\
= 2 \, \sum_{j \in (\Z/8\Z)^*} \chi_D(j) \, e^{j \al_4}  \, . 
\end{multline*}
Putting these contributions together we obtain
\[  \inv_D(e^{\gamma}) = \frac{1}{48} 
\, \Bigg\{
\sum_{\mu \in M(\gamma)^+} e^{\mu} - \sum_{\mu \in M(\gamma)^-} e^{\mu}
\Bigg\}  \, .  \]
This proves the proposition. \eop

\medskip
\noindent
Note that the decomposition $I \backslash I_4 = M(\gamma)^+ \cup M(\gamma)^-$
% into two disjoint subsets
is independent of the choice of $\gamma$ because $\C[D]^{\Gamma}$ is $1$-dimensional. The sets $M(\gamma)^{\pm}$ have size $2\cdot 8 + 2\cdot 4 = 24 = 48/2$.

\medskip

We remark that every discriminant form $D$ of level $8$, exponent $8$, order $8^3$ and even signature $1+t \! \mod 8$ is isomorphic to $2_1^{+1} 4_t^{\epsilon} \, 8_{I\!I}^{+2}$ with $\epsilon = \left( \frac{t}{2} \right)$.

\section{Induction} \label{heidegger}

In this section we recall some properties of the isotropic induction.

\medskip

Let $D$ be a discriminant form of even signature and $H$ an isotropic subgroup of $D$. Then $H^{\perp}/H$ is a discriminant form of the same signature as $D$ and of order $|H^{\perp}/H| = |D|/|H|^2$. There is an isotropic lift
\[ \uparrow_H^D : \C[H^{\perp}/H] \to \C[D]  \]
defined by
\[  \uparrow_H^D(e^{\gamma+H}) = \sum_{\mu \in H} e^{\gamma + \mu}   \]
for $\gamma\in H^{\perp}$ and an isotropic descent 
\[ \downarrow_H^D : \C[D] \to \C[H^{\perp}/H] \]
defined by
\[ \downarrow_H^D(e^{\gamma}) = 
\begin{cases}
\, e^{\gamma+H} & \text{if $\gamma\in H^{\perp}$}, \\
\,          0 & \text{otherwise} 
\end{cases} 
\]
(see \cite{Br}, section 5.1, \cite{S3}, section 4 and  \cite{S4}, section 2). The following results are easy to prove.

\begin{prp}
  Let $D$ be a discriminant form of even signature and $H$ an isotropic subgroup of $D$. Then:
\begin{itemize}
\item[i)]
  $(v,\uparrow_H^D w)_{\C[D]} = (\downarrow_H^D v,w)_{\C[H^{\perp}/H]}$ for all $v \in \C[D], w \in \C[H^{\perp}/H]$.
\item[ii)]
The diagram   
\[
\begin{tikzcd}
  \C[H^{\perp}/H] \arrow[d, "\rho_{H^{\perp}/H}(M)" '] \arrow[r, "\uparrow_H^D"] & \C[D] \arrow[d,"\rho_{D}(M)"] \\
  \C[H^{\perp}/H] \arrow[r, "\uparrow_H^D"] & \C[D]
\end{tikzcd} 
\]
commutes for all $M \in \Gamma$ and analogously for $\downarrow_H^D$. In particular $\uparrow_H^D$ and $\downarrow_H^D$ map invariants to invariants.
\item[iii)]
  The maps $\uparrow_H^D$ and $\downarrow_H^D$ commute with $\inv_{H^{\perp}/H}$ and $\inv_D$.
\end{itemize}
\end{prp}

\noindent
Let $D$ be a discriminant form of even signature and $H$ an isotropic subgroup of $D$. Recall that the isotropic subgroups of $H^{\perp}/H$ are of the form $K/H$ for some isotropic subgroup $K$ of $D$ with $H \subset K \subset K^{\perp} \subset H^{\perp}$. The isotropic lift is transitive in the following sense.

\begin{prp} \label{amk}
Let $D$ be a discriminant form of even signature and $H, K$ isotropic subgroups of $D$ with $H \subset K$. Then 
$H \subset K \subset K^{\perp} \subset H^{\perp}$ and $K/H$ is an isotropic subgroup of $H^{\perp}/H$ with orthogonal complement $K^{\perp}/H$. The quotient $(K^{\perp}/H)\big/(K/H)$ is naturally isomorphic to the discriminant form $K^{\perp}/K$ and the diagram
\[
  \begin{tikzcd}
    \C\big[\frac{K^{\perp}/H}{K/H}\big] \arrow[r,"\uparrow_{K/H}^{H^{\perp}/H}"] & \C[H^{\perp}/H] \arrow[r, "\uparrow_H^D"] & \C[D]  \\
%    \C[K^{\perp}/K] \arrow[u] \arrow[urr,bend right=18,"\uparrow_K^D"'] 
\C[K^{\perp}/K] \arrow[u] \arrow[urr,"\uparrow_K^D"'] 
  \end{tikzcd} 
\]
commutes.
\end{prp}
Identifying $\C[K^{\perp}/K]$ with $\C\big[\frac{K^{\perp}/H}{K/H}\big]$ we can write
%$\uparrow_H^D \, \circ \, \uparrow_{K/H}^{H^{\perp}/H} = \, \uparrow_K^D$.
\[   \uparrow_H^D \, \circ \, \uparrow_{K/H}^{H^{\perp}/H} \quad = \quad \uparrow_K^D   \, . \]

\section{The main theorem}

In this section we prove the main result of this paper. We define fundamental invariants and show that each invariant is induced from these invariants.

\medskip

As explained in the introduction we can restrict to $p$-adic discriminant forms. Let $D$ be a discriminant form of level $p^l$ where $p$ is a prime and even signature. For $\gamma \in D$ we define 
\[ a(p,\gamma)
  = | \{ H \subset \gamma^{\perp} \text{ is an isotropic subgroup of $D$ with $|H| = p$} \} | \, .  \]
% $a(p,\gamma)$ as the number of isotropic subgroups of order $p$ in $\gamma^{\perp}$.

%\medskip

First we give a sufficient condition for an element $e^{\gamma} \in \C[D]$ to be a linear combination of isotropic lifts (cf.\ also \cite{Wr}, Lemma 120 and Lemma 121).

\begin{prp}  \label{nazare1}
% Let $D$ be a discriminant form of level $q=p^l$, $p$ prime and even signature.
Let $\gamma \in D$. Suppose that $\gamma^{\perp}$ contains an isotropic subgroup $H$ isomorphic to $(\Z/p\Z)^2$. Then $e^{\gamma}$ can be written as a linear combination of isotropic lifts.
\end{prp}
{\em Proof:}
The group $H$ has $p^2-1$ elements of order $p$ and therefore has $p+1 = (p^2-1)/(p-1)$ subgroups of order $p$. We denote them by $H_0,H_1,\ldots,H_p$. The inclusions $H_j \subset H \subset \la \gamma \ra^{\perp}$ imply $\la \gamma \ra \subset H^{\perp} \subset H_j^{\perp}$. Define
\[  v 
= \, \sum_{j=1}^p \uparrow_{H_j}^D(e^{\gamma+H_j}) 
= \, \sum_{j=1}^p \, \sum_{\mu \in H_j} e^{\gamma+\mu} 
= p e^{\gamma} + \, \sum_{\mu \in H\backslash H_0} e^{\gamma + \mu}
\]
and
\[ w 
=  \, \uparrow_{H_0}^D \Big( \, \sum_{\mu \in H_1\backslash \{ 0 \} } e^{\gamma + \mu + H_0} \Big)
= \, \sum_{\mu \in H_1\backslash \{ 0 \} } \, \sum_{\bt \in H_0} e^{\gamma + \mu + \bt}  
= \, \sum_{\mu \in H \backslash H_0} e^{\gamma + \mu} \, .
\]
Then $e^{\gamma} = (v-w)/p$. \eop

\begin{prp} \label{nazare2}
  % Let $D$ be a discriminant form of even signature and $p$ a prime. Let $\gamma \in I \backslash \{ 0 \}$ be of order $q$ where $q$ is a power of $p$.
Let $\gamma \in I \backslash \{ 0 \}$. Then $\gamma^{\perp}$ contains an isotropic subgroup isomorphic to $(\Z/p\Z)^2$ if and only if $a(p,\gamma) > 1$.
\end{prp}
{\em Proof:} Suppose $a(p,\gamma) > 1$. Let $\gamma$ be of order $n$. Then $(n/p) \gamma$ generates an isotropic subgroup of order $p$ in $\gamma^{\perp}$. Since $a(p,\gamma) > 1$, there is another isotropic subgroup of order $p$ in $\gamma^{\perp}$. Both groups together generate an isotropic subgroup isomorphic to $(\Z/p\Z)^2$ in $\gamma^{\perp}$. The other direction is clear. \eop

\medskip

Now we have to distinguish between even and odd primes. The following result is well-known and easy to prove.

\begin{prp}
  Let $D$ be a discriminant form of level $p^l$ where $p$ is an odd prime. Suppose $D$ contains no non-trivial isotropic elements. Then $D$ is isomorphic to one of the following discriminant forms:
%  \[   0, \, p^{\pm 1}, \, p^{\epsilon 2} \text{ with } \epsilon = -\Big( \frac{-1}{p} \Big) \, .  \]
  \[   0, \, p^{\pm 1}, \, p^{\epsilon 2} \text{ with } \epsilon = -\left( \frac{-1}{p} \right) .  \]
\end{prp}
Recall that for these discriminant forms $\dim(\C[D]^{\Gamma})=0$ (see the comment after Proposition \ref{iggypop}).

\begin{prp} \label{fdpodd}
Let $D$ be a discriminant form of level $p^l$ where $p$ is an odd prime. Let $\gamma \in I$ be of order $p$. If $a(p,\gamma)=1$, then $\inv_D(e^{\gamma}) = \inv_D(e^0)$ or $D$ is of type $p^{\epsilon 2}$ with $\epsilon = \big( \frac{-1}{p} \big)$, $p^{\pm 3}$ or $p^{-4}$. 
\end{prp}
{\em Proof:} 
First we consider the case $\gamma \notin D^p$. We show that $\la \gamma \ra^{\perp}/\la \gamma \ra$ contains no non-trivial isotropic elements. Suppose $\mu + \la \gamma \ra \in \la \gamma \ra^{\perp}/\la \gamma \ra$ with $\mu \notin \la \gamma \ra$ is isotropic. Then $\mu \in \la \gamma \ra^{\perp}$ is isotropic and $a(p,\gamma) =1$ implies
%$\frac{n}{p}\mu \in \la \gamma \ra$
$(n/p)\mu \in \la \gamma \ra$
where $n$ is the order of $\mu$. Since $\gamma \notin D^p$,  we conclude $\mu \in \la \gamma \ra$. It follows that $\la \gamma \ra^{\perp}/\la \gamma \ra$ is of type $0$, $p^{\pm 1}$ or $p^{\varepsilon 2}$ with $\varepsilon = -\big( \frac{-1}{p} \big)$. If $\la \gamma \ra^{\perp}/\la \gamma \ra = 0$, then $|D| = p^2$ so that $D$ is isomorphic to $q^{\pm 1}$ with $q=p^2$ or to $p^{\epsilon 2}$ with $\epsilon = \big( \frac{-1}{p} \big)$. The first case contradicts $\gamma \notin D^p$. Hence $D$ is isomorphic to $p^{\epsilon 2}$. If $\la \gamma \ra^{\perp}/\la \gamma \ra = p^{\pm 1}$, then $|D| = p^3$ and $D$ must be of type $p^{\pm 3}$. For $\la \gamma \ra^{\perp}/\la \gamma \ra \cong p^{\varepsilon 2}$ we find $D \cong p^{-4}$.

Now let $\gamma\in D^p$. We choose a Jordan decomposition of $D$ and write $D = A \oplus B$ where $A$ denotes the sum over the irreducible components of exponent $p$ and $B \neq 0$ the sum over the remaining components.
%Then $\gamma \in B^p$. Recall that $B^p$ is the orthogonal complement of $B_p$.
Then $\gamma \in B^p$ and $B_p \subset \gamma^{\perp}$. (Recall that $B_p$ is orthogonal to $B^p$.)  
Since $B_p$ is isotropic and $a(p,\gamma) =1$, we have $B_p = \la \gamma \ra$, i.e.\ $B_p$ is cyclic. This implies that $B$ is cyclic. Let $B \cong q^{\pm 1}$. Then $\gamma = (q/p) \bt$ for some generator $\bt$ of $B$. An isotropic element in $D$ is of the form $\mu + m \beta$ with $\mu \in A$ and $p|m$. Since
% $(\gamma, \mu + m\bt) = (q/p)x(\beta,\beta) = 0 \bmod 1$
\[  (\gamma, \mu + m\bt) = (q/p)m(\beta,\beta) = 0 \! \mod 1 \, ,   \]
this implies $\gamma \in I^{\perp}$. Hence $\inv_D(e^{\gamma}) = \inv_D(e^0)$ by Proposition \ref{iggypop}. \eop

\medskip
%\noindent
We leave the proof of the next result to the reader.

\begin{prp}  \label{nickcave}
  Let $D$ be a discriminant form of level $2^l$. Suppose that $D$ contains no non-trivial isotropic elements. Then $D$ is isomorphic to one of the following discriminant forms:
  	\begin{align*}
		&0, \, 2_{I\!I}^{-2}, \, 2_{t}^{\pm 1}, \\
		&2_{t}^{\pm 2} \text{ with } t=2 \, \bmod 4 \, , \\
		&2_{t}^{\epsilon 3} \text{ with } \epsilon \left( \frac{t}{2} \right) = -1 \, , \\[1mm]
		&4_{t}^{\pm 1}, \, 2_{s}^{\pm 1}4_{t}^{\pm 1} \, .
	\end{align*}
      \end{prp}
      For these discriminant forms the subspace of invariants is trivial. We also remark that we do not assume in the proposition that $D$ has even signature.
      
\begin{prp} \label{fdpeven2}
Let $D$ be a discriminant form of level $2^l$ such that $\chi_D$ is trivial. Let $\gamma \in I$ be of order $2$. If $a(2,\gamma) =1$, then $\inv_D(e^{\gamma}) = \inv_D(e^0)$ or $D$ is of type $2_{I\!I}^{+2}$ or $2_{I\!I}^{-4}$. 
\end{prp}
{\em Proof:} 
Note that the condition on $\chi_D$ implies that $|D|$ is a square and $\sign(D) = 0 \! \mod 4$.

First we consider the case $\gamma \notin D^2$. The discriminant form $\la \gamma \ra^{\perp}/\la \gamma \ra$ has the same signature and square class as $D$ and contains no non-trivial isotropic elements. Hence $\la \gamma \ra^{\perp}/\la \gamma \ra$ is isomorphic to $0$ or $2_{I\!I}^{-2}$. If $\la \gamma \ra^{\perp}/\la \gamma \ra \cong 0$, then $|D|= 2^2$ and $D$ contains a non-trivial isotropic element of order $2$. This implies $D \cong 2_{I\!I}^{+2}$ or $D \cong 2_0^{+2}$. In the latter case $\C[D]^{\Gamma}$ is spanned by $\inv_D(e^{\gamma}) = \inv_D(e^0)$ (cf.\ Proposition \ref{pba}). If $\la \gamma \ra^{\perp}/\la \gamma \ra \cong 2_{I\!I}^{-2}$, then $D$ has order $16$ and signature $4 \, \bmod 8$. The discriminant forms of order $16$ and signature $4 \, \bmod 8$ are
\begin{align*}
%	&2_{t}^{\pm1}8_{s}^{\pm1},\ 4_{4}^{-2},\ 4_{I\!I}^{-2}, \\
%	&2_{t}^{\pm4}, \ 2_{I\!I}^{-4}
4_{4}^{-2}, \, 2_{4}^{+4}, \, 2_{I\!I}^{-4} \, . 
\end{align*}
In the first case the isotropic elements are multiples of $2$. In the case $2_{4}^{+4}$ the space $\C[D]^{\Gamma}$ is trivial so that $\inv_D(e^{\gamma}) = \inv_D(e^0) = 0$ (cf.\ Proposition \ref{pba}).

Next we assume that $\gamma\in D^2$. We choose a Jordan decomposition of $D$ and write $D = A \oplus B$ where $A$ denotes the sum over the irreducible components of exponent $2$ and $B \neq 0$ the sum over the remaining components. Then $\gamma \in B^2$. The group $B^2$ is the orthogonal complement of $B_2 \subset B^2$, but in general $B_2$ is not isotropic. If $B_2$ is isotropic, we can argue exactly as in the proof of Proposition \ref{fdpodd}.
Suppose $B_2$ is not isotropic. Since $a(2,\gamma) =1$, the only non-trivial isotropic element in $B_2$ is $\gamma$. Hence the discriminant form $B$ must be of type $4_t^{\pm 2}$ or $4_s^{\pm 1} q_t^{\pm 1}$ with $8|q$.
%In the latter case we can choose a generator $\bt$ of $q_t^{\pm 1}$ such that $\gamma = (q/2) \bt$. Then $\gamma \in I^{\perp}$ so that $\inv_D(e^{\gamma}) = \inv_D(e^0)$ by Proposition \ref{iggypop}.
In the latter case we can choose a generator $\bt$ of $q_t^{\pm 1}$. Then $\gamma = (q/2) \bt$ and $\gamma \in I^{\perp}$ so that $\inv_D(e^{\gamma}) = \inv_D(e^0)$ by Proposition \ref{iggypop}.
Suppose $B$ is of type $4_t^{\pm 2}$. We choose orthogonal generators $\bt_1, \bt_2$ of $B$. Then $\gamma = 2\bt_1 + 2\bt_2$ and any isotropic element in $D$ is of the form $\mu + m_1\bt_1 + m_2\bt_2$ with $\mu \in A$ and $2| (m_1 + m_2)$. Now 
\[ (\gamma,\mu+m_1\beta_1+m_2\beta_2) = 2m_1(\beta_1,\beta_1) + 2m_2(\beta_2,\beta_2) = 0 \! \mod 1, \]
implies $\gamma \in I^{\perp}$ so that again $\inv_D(e^{\gamma}) = \inv_D(e^0)$ by Proposition \ref{iggypop}.
\eop

\begin{prp} \label{fdpeven4}
Let $D$ be a discriminant form of level $2^l$ and even signature such that $\chi_D$ is non-trivial and $|D|$ is a square. Let $\gamma \in I$ be of order $4$. If $a(2,\gamma)=1$, then $\inv_D(e^{\gamma}) = 0$ or $D$ is of type $2_t^{+2}4_{I\!I}^{+2}$ with $t = 2 \! \mod 4$.
\end{prp}
{\em Proof:}
Since $\chi_D$ is non-trivial and $|D|$ is a square, we have $\sign(D) = 2 \, \bmod 4$.

First we consider the case $\gamma \notin D^2$. The discriminant form $\la \gamma \ra^{\perp}/\la \gamma \ra$ has the same signature and square class as $D$ and contains no non-trivial isotropic elements. Hence it is isomorphic to $2_{t}^{\pm2}$ with $t = 2 \, \bmod 4$. It follows that $D$ has order $64$. The discriminant forms of order $64$ and signature $2 \, \bmod 4$ containing elements of order $4$ are
\begin{align*}
	&2_{s}^{\pm 1} 32_{t}^{\pm 1}, \, 4_{s}^{\pm 1}16_{t}^{\pm 1}, \, 8_{t}^{\pm 2}, \\
	&2_{s}^{\pm 3}8_{t}^{\pm 1}, \, 2_{s}^{\pm 2}4_{t}^{\pm 2}, \, 2_{I\!I}^{\pm 2}4_{t}^{\pm 2}, \\
	&2_{s}^{\pm 2}4_{I\!I}^{\pm 2}
\end{align*}
with suitable $s$, $t$ and signs.
% For the discriminant forms of type $2_{s}^{\pm 1} 32_{t}^{\pm 1}$, $4_{s}^{\pm 1}16_{t}^{\pm 1}$ and $8_{t}^{\pm 2}$ any isotropic element of order $4$ is in $D^2$.
For the discriminant forms of type $2_{s}^{\pm 1} 32_{t}^{\pm 1}$, $4_{s}^{\pm 1}16_{t}^{\pm 1}$ and $8_{t}^{\pm 2}$ the isotropic elements of order $4$ are multiples of $2$. For the discriminant forms of type $2_{s}^{\pm 3}8_{t}^{\pm 1}$, $2_{s}^{\pm 2}4_{t}^{\pm 2}$ and $2_{I\!I}^{\pm 2}4_{t}^{\pm 2}$ any isotropic element $\mu$ of order $4$ satisfies $ 2 \mu \in I^{\perp}$ so that $\inv_D(e^{\mu}) = 0$ by Proposition \ref{ironmaiden}. Finally $2_{s}^{\pm 2} 4_{I\!I}^{\pm 2} \cong 2_{t}^{+2} 4_{I\!I}^{+2}$ for some $t$ with $t = 2 \! \mod 4$.

Now suppose $\gamma\in D^2$. As above we choose a Jordan decomposition of $D$ and write $D = A \oplus B$ where $A$ denotes the sum over the irreducible components of exponent dividing $4$ and $B \neq 0$ the sum over the remaining components. Then $\gamma$ is orthogonal to $B_2$. Since $B_2$ is isotropic and $a(2,\gamma) = 1$, we have $B_2 = \la 2 \gamma \ra$. Hence $B$ is cyclic.
%We can choose a generator $\bt$ of $B \cong q_t^{\pm 1}$ such that $\gamma = 2 \al + (q/4)\bt$ for some $\al \in A$.
We choose a generator $\bt$ of $B \cong q_t^{\pm 1}$. Then $\gamma = 2 \al + (q/4)\bt$ for some $\al \in A$.
An isotropic element in $D$ is of the form $\mu+m\beta$ with $\mu\in A$ and $2|m$. Now
\[  (2 \gamma, \mu + m \bt) = (q/2)m(\bt,\bt) = 0 \! \mod1  \]
so that $2 \gamma \in I^{\perp}$. Hence $\inv_D(e^{\gamma}) = 0$ by Proposition \ref{ironmaiden}. \eop

\begin{prp} \label{fdpeven8}
Let $D$ be a discriminant form of level $2^l$ and even signature such that $|D|$ is not a square. Let $\gamma \in I$ be of order $8$. If $a(2,\gamma)=1$, then $\inv_D(e^{\gamma}) = 0$ or $D$ is of type $2_1^{+1} 4_t^{\epsilon} \, 8_{I\!I}^{+2}$ with $t = 1 \, \bmod 2$ and $\epsilon = \left( \frac{t}{2} \right)$. 
\end{prp} 
{\em Proof:}
As before we consider first the case that $\gamma \notin D^2$. The discriminant form $\la \gamma \ra^{\perp}/ \la \gamma \ra$ has the same signature and square class as $D$ and contains no non-trivial isotropic elements. Hence $\la \gamma \ra^{\perp}/ \la \gamma \ra$ is of type $2_{s}^{\pm 1}4_{t}^{\pm1}$. It follows that $D$ has order $512$. The discriminant forms of order $512$ and even signature containing elements of order $8$ are
\begin{align*}
	&2_{s}^{\pm 1} 256_{t}^{\pm 1}, \, 4_{s}^{\pm 1} 128_{t}^{\pm 1}, \, 8_{s}^{\pm 1} 64_{t}^{\pm 1}, \, 2_{s}^{\pm 3} 64_{t}^{\pm 1}, \\
	&16_{s}^{\pm 1} 32_{t}^{\pm 1}, \, 2_{r}^{\pm2} 4_{s}^{\pm 1} 32_{t}^{\pm 1}, \, 2_{I\!I}^{\pm 2} 4_{s}^{\pm 1} 32_{t}^{\pm 1}, \, 2_{r}^{\pm 2} 8_{s}^{\pm 1} 16_{t}^{\pm 1},\\
	&2_{I\!I}^{\pm 2} 8_{s}^{\pm 1} 16_{t}^{\pm 1}, \, 2_{r}^{\pm 1} 4_{s}^{\pm2} 16_{t}^{\pm 1}, \, 2_{s}^{\pm 1}4_{I\!I}^{\pm2} 16_{t}^{\pm 1}, \, 2_{s}^{\pm 5} 16_{t}^{\pm 1}, \\
	&2_{r}^{\pm 1} 4_{s}^{\pm 1} 8_{t}^{\pm 2}, \, 2_{s}^{\pm 1} 4_{t}^{\pm 1} 8_{I\!I}^{\pm 2}, \, 4_{s}^{\pm 3} 8_{t}^{\pm 1}, \, 2_{r}^{\pm 4} 4_{s}^{\pm 1} 8_{t}^{\pm 1}, \\ 
	&2_{I\!I}^{\pm 4} 4_{s}^{\pm 1} 8_{t}^{\pm 1}
\end{align*}
with suitable $s$, $t$ and signs.
In the discriminant forms $2_{s}^{\pm 1} 256_{t}^{\pm 1}$ and $4_{s}^{\pm 1} 128_{t}^{\pm 1}$ the isotropic elements of order $8$ are multiples of $2$ contradicting our assumption on $\gamma$.
%If $D$ is of type $2_{s}^{\pm 1} 256_{t}^{\pm 1}$ or $4_{s}^{\pm 1} 128_{t}^{\pm 1}$ the isotropic elements of order $8$ in $D$ are in $D^2$ contradicting our assumption on $\gamma$.
The discriminant forms
\begin{align*}
	&8_{s}^{\pm 1} 64_{t}^{\pm 1}, \, 16_{s}^{\pm 1} 32_{t}^{\pm 1}, \, 2_{r}^{\pm 2} 8_{s}^{\pm 1} 16_{t}^{\pm 1},\\
	&2_{I\!I}^{\pm 2} 8_{s}^{\pm 1} 16_{t}^{\pm 1}, \, 2_{s}^{\pm 1} 4_{I\!I}^{\pm 2} 16_{t}^{\pm 1}, \, 2_{s}^{\pm 5}16_{t}^{\pm 1}, \\
	&4_{s}^{\pm3} 8_{t}^{\pm 1}, \, 2_{r}^{\pm4} 4_{s}^{\pm 1} 8_{t}^{\pm 1}, \, 2_{I\!I}^{\pm 4} 4_{s}^{\pm 1} 8_{t}^{\pm 1}
\end{align*}
contain no isotropic elements of order $8$ so $D$ cannot be isomorphic to any of them.
If $D$ is of type
\begin{align*}
	&2_{s}^{\pm 3} 64_{t}^{\pm 1}, \, 2_{r}^{\pm 2} 4_{s}^{\pm 1} 32_{t}^{\pm 1}, 2_{I\!I}^{\pm 2} 4_{s}^{\pm 1} 32_{t}^{\pm 1},\\
	&2_{r}^{\pm 1} 4_{s}^{\pm 2} 16_{t}^{\pm 1}, 2_{r}^{\pm 1} 4_{s}^{\pm 1} 8_{t}^{\pm 2},
\end{align*}
%\[    2_{s}^{\pm 3} 64_{t}^{\pm 1}, \, 2_{r}^{\pm 2} 4_{s}^{\pm 1} 32_{t}^{\pm 1},
%      2_{I\!I}^{\pm 2} 4_{s}^{\pm 1} 32_{t}^{\pm 1}, 2_{r}^{\pm 1} 4_{s}^{\pm 2} 16_{t}^{\pm 1}
%      \text{ or } 2_{r}^{\pm 1} 4_{s}^{\pm 1} 8_{t}^{\pm 2},    \]
then any isotropic element $\mu$ of order $8$ in $D$ satisfies $ 4 \mu \in I^{\perp}$ so that $\inv_D(e^{\mu}) = 0$ by Proposition \ref{cash}. Finally $2_{r}^{\pm 1} 4_{s}^{\pm 1} 8_{I\!I}^{\pm 2} \cong2_{1}^{+1}4_{t}^{\epsilon} 8_{I\!I}^{+2}$ for some $t$ with $t = 1 \! \mod 2$ and $\epsilon = \left(\frac{t}{2}\right)$.

Now suppose $\gamma\in D^2$. Again we choose a Jordan decomposition of $D$ and write $D = A \oplus B$ where $A$ denotes the sum over the irreducible components of exponent dividing $8$ and $B \neq 0$ the sum over the remaining components. Then $\gamma$ is orthogonal to $B_2$.
Since $B_2$ is isotropic, we have $B_2 = \la 4\gamma \ra$. Hence $B$ is cyclic.
We choose a generator $\bt$ of $B \cong q_t^{\pm 1}$. Then $\gamma = 2 \al + (q/8)\bt$ for some $\al \in A$. An isotropic element in $D$ is of the form $\mu + m\bt$ with $\mu\in A$ and $2|m$. Since 
\[  (4 \gamma, \mu + m \bt) = (q/2)m(\bt,\bt) = 0 \! \mod 1 \, ,  \]
this implies $4 \gamma \in I^{\perp}$. Hence $\inv_D(e^{\gamma}) = 0$ by Proposition \ref{cash}.  \eop

\medskip

The above discriminant forms, with the exception of $p^{\epsilon 2}$ and $2_{I\!I}^{+2}$, play an important role in our main result. We summarise some of their properties in the following tables. First let $p$ be an odd prime: 

\[
\renewcommand{\arraystretch}{1.2}
\begin{array}{c|c|c|c}
 D & \text{square class} & \text{signature} & \text{invariant} \\[0.5mm] \hline 
   &   &   & \\[-4mm]
 0           & \text{square}           & 0 \, \bmod 8      & e^0 \\
 p^{-4}       & \text{square}           & 4 \, \bmod 8      
                       & (p-1) e^0 - \sum_{\gamma \in M} e^{\gamma} \\
 p^{\epsilon 3} & \text{non-square}       & 0 \, \bmod 2      
                       & \sum_{\gamma \in M^+}e^{\gamma} - \sum_{\gamma \in M^-}e^{\gamma} 
\end{array}
\]

\vspace*{2mm}
\noindent
The case $p=2$ is more complicated:

\vspace*{-3mm}
\[
\renewcommand{\arraystretch}{1.2}
\begin{array}{c|c|c|c}
 D & \text{square class} & \text{signature} & \text{invariant} \\[0.5mm] \hline 
   &   &   & \\[-4mm]
 0           & \text{square}           & 0 \, \bmod 8      & e^0 \\
 2_{I\!I}^{-4}       & \text{square}           & 4 \, \bmod 8      
                       &  e^0 - \sum_{\gamma \in M} e^{\gamma} \\
 2_t^{+2}4_{I\!I}^{+2}  &  \text{square}  & t = 2 \, \bmod 4 
                       & \sum_{\gamma \in M^+}e^{\gamma} - \sum_{\gamma \in M^-}e^{\gamma} \\
 2_1^{+1} 4_t^{\epsilon} \, 8_{I\!I}^{+2}  &  \text{non-square}  &  1+t = 0 \, \bmod 2 
                       & \sum_{\gamma \in M^+}e^{\gamma} - \sum_{\gamma \in M^-}e^{\gamma} 
\end{array}
\]
In all these cases $\C[D]^{\Gamma}$ is $1$-dimensional. We wrote $M$ for the set of isotropic elements whose order is equal to the level of $D$. In the indicated cases $M$ has a canonical decomposition $M = M^+ \cup M^-$. We denote the above discriminant forms as 
%$F^{x,s}_p$ 
$D^{x,s}_p$
where $x$ is the square class and $s$ the signature of $D$ and the generator of the subspace of invariants as 
%$f^{x,s}_p$
$i^{x,s}_p$.

\begin{thm} \label{mainth}
Let $D$ be a discriminant form of even signature $s$, square class $x$ and level $p^l$ where $p$ is a prime. Then the invariants of the Weil representation on $\C[D]$ are generated by the invariants $\uparrow_H^D(i^{x,s}_p)$ where $H$ is an isotropic subgroup of $D$ such that $H^{\perp}/H$ is isomorphic to the discriminant form $D^{x,s}_p$.
\end{thm}
{\em Proof:}
Recall that the invariants $\inv_D(e^{\gamma})$, $\gamma \in I$ generate $\C[D]^{\Gamma}$. Let $\gamma \in I$. We will show below that at least one of the following statements applies:
\begin{itemize}
\item[i)]
$D$ is a fundamental discriminant form,  
\item[ii)]
$\inv_D(e^{\gamma})$ is induced from smaller discriminant forms of the same signature and square class as $D$, i.e.\ $\inv_D(e^{\gamma})$ is a linear combination of lifts of invariants for suitable isotropic subgroups of $D$,
\item[iii)]
$\inv_D(e^{\gamma})=0$.
\end{itemize}
Then the theorem follows by induction on the order of $D$: If $|D| = 1$, the discriminant form is fundamental. Let $|D|>1$. If $D$ is fundamental, there is nothing to prove. Suppose $D$ is not fundamental. Let $\gamma\in I$. If $\inv_D(e^{\gamma}) \neq 0$, then
%$\inv_D(e^{\gamma})$
it is a linear combination of invariants which are lifts of invariants on smaller discriminant forms. The induction hypothesis implies that the invariants on the smaller discriminant forms are induced from the fundamental invariant corresponding to $D$. By the transitivity of the isotropic lift (see Proposition \ref{amk}) $\inv_D(e^{\gamma})$ is a linear combination of lifts of the fundamental invariant on isotropic subgroups of $D$. This finishes the induction.

Now we prove that at least one of the above three statements holds. We assume that $D$ is non-trivial. If $D$ contains no non-trivial isotropic elements, then $\inv_D(e^{\gamma}) = 0$ for all $\gamma \in D$.
% We now assume that $I \neq \{ 0 \}$.
Suppose $I \neq \{ 0 \}$. Let $\gamma \in I$.
Since 
%Then
%\[  \rho_D(S) e^0  = \frac{e(\sign(D)/8)}{\sqrt{|D|}} \sum_{\beta \in D} e^{\beta}   \]
%so that
\[  \inv_D(e^0) = \rho_D(S) \inv_D ( e^0 ) = \inv_D (\rho_D(S) e^0 ) = \frac{e(\sign(D)/8)}{\sqrt{|D|}} 
 \sum_{\beta \in I} \inv_D(e^{\beta})  \, ,  \] 
we can assume that $\gamma \neq 0$.
 
We define $m=p$ if $p$ is odd and 
\[ m = 
\begin{cases}
\, 2 & \text{if $|D|$ is a square and $\sign(D) = 0 \, \bmod 4$,} \\
\, 4 & \text{if $|D|$ is a square and $\sign(D) = 2 \, \bmod 4$,} \\
\, 8 & \text{if $|D|$ is a non-square} \\       
\end{cases} 
\]
for $p=2$.

First we consider the case that $\gamma \notin D_m$. Let $n$ be the order of $\gamma$ and $H = \la \gamma \ra_p$. Then for all $v = \sum_{\bt \in D} v_{\bt} e^{\bt} \in \C[D]^{\Gamma}$ we have
\begin{gather*}
(v, \uparrow_H^D( \inv_{H^{\perp}/H}(e^{\gamma + H}) ) 
= (v, \inv_D( \uparrow_H^D(e^{\gamma + H}) ) )
= \sum_{\mu \in H} (v, \inv_D( e^{\gamma + \mu}) ) \\
= \sum_{\mu \in H} v_{\gamma + \mu} 
= \sum_{\substack{ a \in \Z/n\Z \\ a = 1 \bmod n/p}} v_{a\gamma}
= \sum_{\substack{ a \in \Z/n\Z \\ a = 1 \bmod n/p}} \chi_D(a) v_{\gamma} 
= p v_{\gamma}
= p (v,\inv_D(e^{\gamma}))
\end{gather*}
because $m|\frac{n}{p}$ so that 
\[   \inv_D(e^{\gamma}) 
= \frac{1}{p} \uparrow_H^D( \inv_{H^{\perp}/H}(e^{\gamma + H})) \, .  \]

Next we consider the case $\gamma \in D_m\backslash \{0 \}$. If $e^{\gamma}$ is a linear combination of isotropic lifts for suitable isotropic subgroups, then the same holds for $\inv_D(e^{\gamma})$ because isotropic induction and $\inv$ commute. We assume that $e^{\gamma}$ is not a linear combination of isotropic lifts. Then $a(p,\gamma) = 1$ by Propositions \ref{nazare1} and \ref{nazare2}.

Suppose $\chi_D$ is trivial. Then $m=p$ and $\inv_D(e^{\gamma})= \inv_D(e^0)$ or $D$ is of type $p^{\epsilon 2}$ with $\epsilon =\big( \frac{-1}{p} \big)$ or $p^{-4}$ if $p$ is odd or of type $2_{I\!I}^{+2}$ or $2_{I\!I}^{-4}$ if $p=2$ (see Propositions \ref{fdpodd} and \ref{fdpeven2}). We go through the possible cases. If $\inv_D(e^{\gamma})= \inv_D(e^0)$ define $H = \la \gamma \ra$. Then for all $v = \sum_{\bt \in D} v_{\bt} e^{\bt} \in \C[D]^{\Gamma}$ we have
\begin{gather*}
(v, \uparrow_H^D( \inv_{H^{\perp}/H}(e^{0 + H}) ) 
  = \sum_{\bt \in H} v_{\bt}
  = v_0 + \sum_{a \in (\Z/p\Z)^*} v_{a\gamma} 
  = v_0 + (p-1) v_{\gamma} \\
  = (v,\inv_D(e^0)) + (p-1) (v,\inv_D(e^{\gamma}))
  = p (v,\inv_D(e^{\gamma}))
\end{gather*}
by Proposition \ref{abba} so that 
$\inv_D(e^{\gamma}) = \frac{1}{p} \uparrow_H^D( \inv_{H^{\perp}/H}(e^{0 + H}) )$.
If $D$ is of type $p^{\epsilon 2}$ with $\epsilon =\big( \frac{-1}{p} \big)$, then $\C[D]^{\Gamma}$ is generated by the characteristic functions of the $2$ maximal isotropic subgroups (see the example after Theorem \ref{dimformulaodd}). The same analysis holds for $D$ of type $2_{I\!I}^{+2}$. The cases $p^{-4}$ and $2_{I\!I}^{-4}$ correspond to fundamental discriminant forms.

Finally we assume that $\chi_D$ is non-trivial.
If $m=p$ is odd, then $\inv_D(e^{\gamma})= \inv_D(e^0) = 0$ or $D$ is of type $p^{\pm 3}$ (see Propositions \ref{tamtam} and \ref{fdpodd}).
Suppose $m=4$. Then $\sign(D) = 2 \! \mod 4$.
If $2\gamma =0$, then $\inv_D(e^{\gamma})=0$ by Proposition \ref{ironmaiden}.
% $v_{\gamma} = \chi_D(3) v_{3\gamma} = -v_{\gamma} = 0$ for all $v = \sum_{\bt \in D} v_{\bt} e^{\bt} \in \C[D]^{\Gamma}$ which implies $\inv_D(e^{\gamma})=0$.
If $\gamma$ has order $4$, then $\inv_D(e^{\gamma})=0$ or $D$ is of type $2_t^{+2}4_{I\!I}^{+2}$ (see Proposition \ref{fdpeven4}).
The case $m=8$ is analogous and uses Proposition \ref{fdpeven8}. \eop

\medskip

A few comments are in order. It is possible that more than one of the conditions i), ii) and iii) applies (see e.g.\ Proposition \ref{vi4}). A consequence of the theorem is that the invariants are defined over $\Z$. A more direct proof of this fact is given in \cite{ES}. The theorem extends Theorem 4.11 in \cite{M} to $p=2$.

\medskip

We describe some examples. Let $p$ be an odd prime and $D$ a discriminant form of even signature. If $|D|=p$, then $D$ is not fundamental and $\dim(\C[D]^{\Gamma}) = 0$.
%The only isotropic element in $D$ is $0$ and $\inv_D(e^0)=0$ by Proposition \ref{tamtam} or Theorem \ref{invelementaryodd}. It follows $\dim(\C[D]^{\Gamma})=0$.
Suppose $|D|=p^2$. Then $D$ is not fundamental and there are three possibilities. If $D$ has level $p$ and is anisotropic, then $\dim(\C[D]^{\Gamma}) = 0$. If $D$ has level $p$ and is isotropic, then $D$ has two non-trivial isotropic subgroups $H_i$, $i=1,2$ of order $p$ with fundamental quotients $H_i^{\perp}/H_i \cong 0$. They generate $\C[D]^{\Gamma}$ which has dimension $2$. If $D$ has level $p^2$, then $D$ has a unique non-trivial isotropic subgroup $H$ with fundamental quotient $H^{\perp}/H \cong 0$. It follows $\dim(\C[D]^{\Gamma}) = 1$. Finally let $|D|=p^3$. We only consider the case that $D$ has level $p$. Then $D$ is fundamental. Nevertheless $D$ has non-trivial isotropic subgroups $H_i$ of order $p$. Here the quotients $H_i^{\perp}/H_i$ have order $p$ so that no non-trivial invariants can be induced from them.

\begin{cor}
Let $D$ be a discriminant form of even signature $s$, square class $x$ and level $p^l$ where $p$ is a prime. Suppose $|D| < |D^{x,s}_p|$. Then $\dim (\C[D]^{\Gamma}) = 0$.
\end{cor}
    {\em Proof:}
If there were non-trivial invariants in $\C[D]$, they would be induced from $D^{x,s}_p$. This is impossible. \eop     
%The space $\C[D]^{\Gamma}$ is generated by the elements $\inv_D(e^{\gamma})$, $\gamma \in I$. Since $D$ is not fundamental, $\inv_D(e^{\gamma})=0$ or $\inv_D(e^{\gamma})$ is induced from $D^{x,s}_p$ (see the proof of Theorem \ref{mainth}). The latter is impossible. \eop

\section{Applications}

The above results have several applications.
For example, the dimension of the space of weight-2 cusp forms transforming under the Weil representation has contributions coming from the invariants.
Furthermore, the theta expansion gives an isomorphism between modular forms for the Weil representation and Jacobi forms of lattice index. The invariants of the Weil representation can be applied to give simple generating sets for Jacobi forms of singular weight.
Another example comes from orthogonal modular forms. Borcherds' additive theta lift (Theorem 14.3 in \cite{B1}) maps the invariants of the Weil representation to orthogonal modular forms of singular weight. This allows to study orthogonal modular forms of singular weight with a special boundary behaviour.
%special behaviour at the 1-dimensional cusps. 
We will describe the first two examples in more detail.

\subsection*{A dimension formula for cusp forms of weight $2$}

Let $D$ be a discriminant form of level $p$ where $p$ is a prime. We give an explicit formula for the dimension of the space $S_2(D)$ of cusp forms of weight $2$ for the Weil representation $\rho_D$.

\medskip

Let $\rho$ be a finite-dimensional representation of $\SL_2(\Z)$ with finite image. Then the dimension of the space of modular forms for $\rho$ of weight at least $2$ can be determined by means of the Selberg trace formula or the Riemann-Roch theorem (see e.g.\  \cite{Sk1}, \cite{B2} and \cite{F}). In weight $2$ there is a contribution coming from the invariants of $\rho$.
%Here we will follow Freitag's approach.
We follow Freitag's approach \cite{F} here.

\medskip

Let $D$ be a discriminant form of prime level. We assume that $D$ is of type $p^{\epsilon n}$ with $n$ even. The argument for odd $n$ is similar. Then $\sign(D) = 0 \, \bmod 4$ so that $Z$ acts as $\rho_D(Z)e^{\gamma} = e^{-\gamma}$. The space $V \subset \C[D]$ spanned by the elements $e^\gamma + e^{-\gamma}$, $\gamma \in D$ is invariant under $\rho_D$. Let $\rho$ be the restriction of $\rho_D$ to $V$ and $d = \dim(V)$. For a complex $d \times d$-matrix $M$ of finite order with eigenvalues $e(x_i)$, $0 \leq x_i < 1$ define
\[  \al(M) = \sum_{i=1}^d x_i \, , \]
in particular
\[ \al(M) =
  \begin{cases}
  \, {\displaystyle \frac{d}{4} - \frac{\tr(M)}{4}} & \text{if $M^2 = I$}, \\[3.5mm]
  \, {\displaystyle \frac{d}{3} - \frac{1}{3}\Rep(\tr(M^{-1})) + \frac{1}{3\sqrt{3}}\Imp(\tr(M^{-1})) }
                                   & \text{if $M^3=I$}.
  \end{cases}
\]
Then the dimension of $S_2(D)$ is given by 
\begin{multline*}
\dim S_2(D) = 
\frac{d}{6} + d - \alpha\big( e(1/2)\rho(S) \big) - \alpha\big( (e(1/3)\rho(ST))^{-1} \big) - \alpha\big( \rho(T) \big)  \\
- | \{ \gamma \in D/\{\pm 1\} | \q(\gamma) = 0 \, \bmod 1 \}| + \dim \C[D]^{\Gamma}  
\end{multline*} 
(see Theorem 6.1 in \cite{F}). We can evaluate this expression using Theorem \ref{dimformulaodd}.

\begin{thm}
Let $D$ be a discriminant form of prime level $p$ and type $p^{\epsilon n}$ with $n$ even. Then $\dim S_2(D) = 0$ if $p \leq 3$ and
\[  \dim S_2(D) = \frac{p^n+5}{24} - \frac{p^{n-1}}{4} - \epsilon \left(\frac{-1}{p} \right)^{n/2} \, \frac{p-5}{4} p^{(n-2)/2} 
%    \epsilon {\displaystyle \left( \frac{-1}{p} \right)^{n/2} } (p-5)p^{(n-2)/2}}{4} 
         + \frac{p^{n-1}-p}{p^2-1}   \]
if $p >3$.
\end{thm}
{\em Proof:} Since $\Gamma(p)$ acts trivial in the Weil representation $\rho_D$, the components of an element in $S_2(D)$ are cusp forms for $\Gamma(p)$. The spaces $S_2(\Gamma(p))$ are trivial for $p \leq 3$ so that $\dim S_2(D) = 0$ in these cases. Suppose $p>3$. Clearly
\[   d = \frac{p^n-1}{2} + 1 = \frac{p^n+1}{2}  \, . \]
Proposition \ref{numberofelementspodd} implies
\begin{align*}
  | \{ \gamma \in D/\{\pm 1\} \, | \, \q(\gamma) = 0 \, \bmod 1 \}|
        & = \frac{N(p^{\epsilon n},0)-1}{2} + 1 \\
	&= \frac{p^{n-1}+1}{2} + \epsilon \left( \frac{-1}{p} \right)^{n/2} \frac{p-1}{2}p^{(n-2)/2} 
\end{align*}
and
\begin{align*}
  \alpha(\rho(T)) &= \sum_{j=0}^{p-1} \, \frac{j}{p} \, | \{ \gamma \in D/\{\pm 1\} \, | \, \q(\gamma) = - j/p \, \bmod 1 \}| \\
  	&= \frac{1}{2} \sum_{j=1}^{p-1} \, \frac{j}{p} \, N(p^{\epsilon n}, - j/p) \\
	&= \frac{p-1}{4} \big( \, p^{n-1} - \epsilon \left( \frac{-1}{p} \right)^{n/2} p^{(n-2)/2} \big) \, . 
\end{align*}
Since $e(1/2)\rho(S)$ has order $2$, we can apply the above formula to calculate $\al(e(1/2)\rho(S))$. For the trace of $e(1/2)\rho(S)$ we find
\begin{align*}
  \tr(e(1/2)\rho(S))
  &= - \frac{1}{4} \sum_{\gamma\in D}(\rho(S)(e^\gamma+e^{-\gamma}),e^\gamma+e^{-\gamma}) \\
  &= - \frac{e(\sign(D)/8))}{4 p^{n/2}} \sum_{ \beta, \gamma \in D} e((\beta,\gamma)) (e^{\beta}+e^{-\beta}, e^{\gamma}+e^{-\gamma}) \\
  &= - \frac{e(\sign(D)/8))}{2 p^{n/2}} \sum_{\gamma \in D} \{ e(2\q(\gamma)) + e(-2\q(\gamma)) \} \\
  &= - \frac{e(\sign(D)/8))}{p^{n/2}} \, e(\sign(D)/8)) p^{n/2} \\
  &= - 1
\end{align*}
where we used Theorem 3.9 in \cite{S2} to evaluate the last sum. Hence
\[  \al\big( e(1/2)\rho(S) \big) = \frac{d}{4} - \frac{\tr(e(1/2)\rho(S))}{4} = \frac{p^n + 3}{8}  \, . \]
Similarly we find
\[  \alpha\big( (e(1/3)\rho(ST))^{-1} \big) = \frac{p^n + 3}{6}  \, . \]
Finally $\dim \C[D]^{\Gamma}$ is given in Theorem \ref{dimformulaodd}. Combining the contributions we obtain the desired formula for the dimension of $S_2(D)$. \eop

\subsection*{Jacobi forms of singular weight}

The space of Jacobi forms $J_{k,L}$ of lattice index $L$ and singular weight $k=\rk(L)/2$
% where $L$ is a positive definite even lattice
is naturally isomorphic to the space of invariants $\C[L'/L]^{\Mp_2(\Z)}$. This allows us to write down a simple generating set for this space.

\medskip

Jacobi forms of lattice index are natural generalizations of Jacobi forms in one variable \cite{EZ}. They were introduced by Gritsenko \cite{G}. Classical examples are Jacobi theta functions. We recall the definition of Jacobi forms of lattice index and describe some of their properties (cf.\ e.g.\ \cite{Sk2}, \cite{GSZ}).

\medskip

The metaplectic group $\Mp_2(\R)$ is the unique connected double cover of the group $\SL_2(\R)$. Its elements can be written as pairs
% $(M,\omega(\tau))$
$(M,\omega)$
where $M = \left(\begin{smallmatrix} a & b \\ c & d \end{smallmatrix}\right)$ is in $\SL_2(\R)$ and $\omega$ is a holomorphic function on the upper half plane ${\cal H}$ such that $\omega(\tau)^2 = c\tau+d$. Then the product of two elements in $\Mp_2(\R)$ is given by
\[   (M_1,\omega_1(\tau)) \, (M_2,\omega_2(\tau))
  = ( M_1M_2,\omega_1(M_2\tau)\omega_2(\tau))  \, .  \]
The inverse image of $\SL_2(\Z)$ under the covering $\Mp_2(\R) \to \SL_2(\R)$ is denoted by $\Mp_2(\Z)$.

Let $L$ be a positive-definite even lattice of rank $n$. Then $\Mp_2(\Z)$ acts from the right on the pairs $(\lambda,\mu) \in L \times L$. The corresponding semidirect product $J_L = \Mp_2(\Z) \ltimes (L \times L)$ is the Jacobi group of lattice index $L$. Recall that the product of two elements in $J_L$ is given by
\begin{align*}
  \MoveEqLeft
  \big( (M_1,\omega_1(\tau)),(\lambda_1,\mu_1) \big) \,
  \big( (M_2,\omega_2(\tau)),(\lambda_2,\mu_2) \big) \\
	& = \big( (M_1M_2,\omega_1(M_2\tau)\omega_2(\tau)),(\lambda_1,\mu_1)M_2 + (\lambda_2,\mu_2) \big) \, .
\end{align*}
We identify $\Mp_2(\Z)$ with the subgroup $\Mp_2(\Z) \ltimes (0 \times 0)$ and write $[\lambda,\mu]$ for the element $(1, (\lambda, \mu)) \in J_L$.

Now let $k \in \frac{1}{2}\Z$. We define an action of the Jacobi group $J_L$ on the functions
$\phi: {\cal H} \times (L\otimes_\Z\C)\rightarrow\C$ by
\begin{align*}
  \phi|_k(M,\omega)(\tau,z)
  &= \phi \bigg(M\tau,\frac{z}{c\tau+d}\bigg) \omega(\tau)^{-2k} e \bigg( \frac{-cz^2/2}{c\tau+d} \bigg) \\
  \phi|_k[\lambda,\mu](\tau,z)
  &= \phi\left(\tau,z + \lambda\tau+\mu\right) e\left(\tau \lambda^2/2 + (\lambda,z)\right)
\end{align*}
where $M = \left(\begin{smallmatrix} a & b \\ c & d \end{smallmatrix}\right) \in \SL_2(\Z)$, $\lambda,\mu \in L$. A Jacobi form of weight $k$ and index $L$ is a holomorphic function $\phi : {\cal H} \times (L\otimes_\Z\C) \to \C$ which is invariant under the action of $J_L$ and possesses a Fourier expansion of the form
\[ \phi(\tau,z)
  = \sum_{\substack{ m \in \Z ,\, \al \in L' \\ m \geq \al^2/2} } c(m,\al) e(m \tau + (\al,z)) \, .  \]
We denote the space of Jacobi forms of weight $k$ and lattice index $L$ by $J_{k,L}$. A Jacobi form $\phi \in J_{k,L}$ has a unique theta expansion
\[  \phi(z,\tau) = \sum_{\gamma \in L'/L} \vartheta_{\gamma}(z,\tau) f_{\gamma}(\tau)  \]
where 
\[  \vartheta_{\gamma}(\tau,z) = \sum_{\al \in \gamma + L} e( \tau \al^2/2 + (\al,z) )  \]
%\[  \vartheta_{\gamma}(\tau,z) = \sum_{x \in \gamma +L} e( \tau \q(\gamma) + (\al,z))  \]
is the Jacobi theta function of the coset $\gamma + L$ and $f(\tau) = \sum_{\gamma \in L'/L} f_{\gamma}(\tau)e^{\gamma}$ a modular form for the Weil representation of $L'/L$ (holomorphic on ${\cal H}$ and at the cusp $\infty$). 
%This decomposition gives a map
We obtain a map
\[  J_{k,L} \to M_{k-n/2}(L'/L)  \]
which is actually an isomorphism. This implies that $J_{k,L}$ is trivial for $k<n/2$. The weight $k=n/2$ is called singular weight. In this case we have an isomorphism
\begin{align*}
   \C[L'/L]^{\Mp_2(\Z)} & \xrightarrow{\varphi_L}     J_{n/2,L} \\
  %  \phi_L : \C[L'/L]^{\Mp_2(\Z)} & \longrightarrow     J_{n/2,L} \\
  \sum_{\gamma \in L'/L} v_{\gamma} e^{\gamma} 
                      & \longmapsto \sum_{\gamma \in L'/L} v_{\gamma} \vartheta_{\gamma}
\end{align*}
(cf.\ also Theorem 5 in \cite{Sk2}). Hence $J_{n/2,L}$ is trivial for odd $n$. For even $n$ we can
describe a simple generating set.
%give a simple generating set.
%Applying Theorem \ref{mainth} we obtain:

\begin{thm}  \label{Iggypop}
  Let $L$ be a positive-definite even lattice of rank $n$ and level $N$. Suppose $n$ is even. For $p|N$ we denote the square class and the signature of the $p$-adic component of $L'/L$ by $x_p$ resp.\ $s_p$. Let $\cal{L}$ be the set of all overlattices $M \supset L$ such that the $p$-adic component of $M'/M$ is isomorphic to $D_p^{x_p,s_p}$ for all $p|N$. Then 
  \[ J_{n/2,L}
    = \sum_{M\in\mathcal{L}} \C \Bigg( \, \sum_{\gamma\in M'/M} v_\gamma\vartheta_{M,\gamma} \bigg)  \]
  where $\sum_{\gamma\in M'/M} v_{\gamma} e^\gamma \in \C[M'/M]^{\SL_2(\Z)}$ is the invariant corresponding to the product
  $\prod_{p|N} i_p^{x_p,s_p}$.
\end{thm}
{\em Proof:}
For $M \in {\cal L}$ we have $L \subset M \subset M' \subset L'$ and $M/L$ is an isotropic subgroup of $L'/L$. Let $v = \sum_{\gamma \in M'/M} v_{\gamma} e^{\gamma} \in \C[M'/M]^{\SL_2(\Z)}$. Then
\[  \uparrow_{M/L}^{L'/L}(v)
  = \sum_{\gamma \in M'/M} v_{\gamma} \uparrow_{M/L}^{L'/L}(e^{\gamma})
  = \sum_{\gamma \in M'/M} v_{\gamma} \sum_{\bt \in M/L} e^{\gamma + \bt}
  \]
  so that
\begin{align*}
  \varphi_L(\uparrow_{M/L}^{L'/L}(v))
  &= \sum_{\gamma\in M'/M} v_{\gamma} \sum_{\beta\in M/L} \vartheta_{L,\gamma+\beta} \\
  &= \sum_{\gamma\in M'/M} v_{\gamma} \sum_{\bt \in M/L}
                  \,  \sum_{\al \in \gamma + \bt + L} e( \tau \al^2/2 + (\al,z) ) \\
  &= \sum_{\gamma\in M'/M} v_{\gamma} \sum_{\al \in \gamma + M} e( \tau \al^2/2 + (\al,z) ) \\
  &= \sum_{\gamma\in M'/M} v_{\gamma} \vartheta_{M,\gamma} \\
  &= \varphi_M(v) 
\end{align*}
because
%$\bigcup_{\beta \in M/L}(\gamma+\beta+L) = \gamma + M$.
\[ \gamma + \bigcup_{\beta\in M/L}( \beta + L ) = \gamma + M \,  .  \]
Hence the diagram 
\[ \begin{tikzcd}
	\C[L'/L]^{\SL_2(\Z)} \arrow{r}{\varphi_L} & J_{n/2,L} \\
	\C[M'/M]^{\SL_2(\Z)} \arrow{r}{\varphi_M} \arrow[swap]{u}{\uparrow_{M/L}^{L'/L}} & J_{n/2,M} \arrow[u,hook]
\end{tikzcd}
\]
commutes. The assertion now follows from Theorem \ref{mainth}.
\eop


\begin{thebibliography}{mmmm}
\addcontentsline{toc}{section}{References}
\bibitem[AGM]{AGM} D.\ Allcock, I.\ Gal, A.\ Mark, {\em The Conway-Sloane calculus for $2$-adic lattices}, Enseign.\ Math.\ {\bf 66} (2020), 5--31 
\bibitem[BEW]{BEW} B.\ C.\ Berndt, R.\ J.\ Evans, K.\ S.\ Williams, {\em Gauss and Jacobi sums}, Canadian Mathematical Society Series of Monographs and Advanced Texts, John Wiley \& Sons, New York, 1998
\bibitem[Bi]{Bi} P.\ Bieker, {\em Invariants for the Weil representation and modular units for orthogonal groups of signature $(2,2)$}, J.\ Number Theory {\bf 250} (2023), 155--182
\bibitem[Bo1]{B1} R.\ E.\ Borcherds, {\em Automorphic forms with singularities on Grassmannians}, Invent.\ Math.\ {\bf 132} (1998), 491--562
\bibitem[Bo2]{B2} R.\ E.\ Borcherds, {\em Reflection groups of Lorentzian lattices}, Duke Math.\ J.\ {\bf 104} (2000), 319--366 
\bibitem[Br]{Br} J.\ H.\ Bruinier, {\em Borcherds products on $\mathrm{O}(2,l)$ and Chern classes of Heegner divisors}, Lecture Notes in Mathematics {\bf 1780}, Springer, Berlin, 2002
\bibitem[CS]{CS} J.\ H.\ Conway, N.\ J.\ A.\ Sloane, {\em Sphere packings, lattices and groups}, 3rd ed., Grundlehren der mathematischen Wissenschaften {\bf 290}, Springer, New York, 1998
\bibitem[DS]{DS} F.\ Diamond, J.\ Shurman, {\em A first course in modular forms}, Graduate Texts in Mathematics {\bf 228}, Springer, New York, 2005
\bibitem[ES]{ES} S.\ Ehlen, N.-P.\ Skoruppa, {\em Computing invariants of the Weil representation}, L-functions and automorphic forms, 81--96, Contrib.\ Math.\ Comput.\ Sci.\ {\bf 10}, Springer, Cham, 2017
\bibitem[EZ]{EZ} M.\ Eichler, D.\ Zagier, {\em The theory of Jacobi forms}, Progress in Mathematics {\bf 55}, Birkh\"auser Boston, Boston, MA, 1985
\bibitem[F]{F} E.\ Freitag, {\em Dimension formulae for vector valued automorphic forms}, Preprint, 2012, available at https:/$\!$/www.mathi.uni-heidelberg.de/$\sim$freitag/papers.html
\bibitem[G]{G} V.\ Gritsenko, {\em Fourier-Jacobi functions in $n$ variables}, J.\ Soviet Math.\ {\bf 53} (1991), 243--252 
\bibitem[GSZ]{GSZ} V.\ Gritsenko, N.-P.\ Skoruppa, D.\ Zagier, {\em Theta blocks}, arXiv: 1907.00188, 2019
\bibitem[M]{M} M.\ K.-H.\ M\"{u}ller, {\em Invariants of the Weil representation of $\SL_2(\Z)$ on discriminant forms of odd order}, Master Thesis,
 Technische Universit\"at Darmstadt, 2021
%Technical University of Darmstadt, 2021
%\bibitem[MS]{MS} S.\ M\"oller, N.\ R.\ Scheithauer, {\em Dimension formulae and generalised deep holes of the Leech lattice vertex operator algebra}, to appear in Ann.\ of Math.
%(arXiv:1910.04947 [math.QA]).
\bibitem[NRS]{NRS} G.\ Nebe, E.\ M.\ Rains, N.\ J.\ A.\ Sloane, {\em Self-dual codes and invariant theory}, Algorithms and Computation in Mathematics {\bf 17},
  % Springer-Verlag, Berlin, 2006
  Springer, Berlin, 2006
\bibitem[N]{N} V.\ V.\ Nikulin, {\em Integral symmetric bilinear forms and some of their applications}, Math.\ USSR Izv.\ {\bf 14} (1980), 103--167
\bibitem[S1]{S1} N.\ R.\ Scheithauer, {\em On the classification of automorphic products and generalized Kac-Moody algebras}, Invent.\ Math.\ {\bf 164} (2006), 641--678
\bibitem[S2]{S2} N.\ R.\ Scheithauer, {\em The Weil representation of $\SL_2(\Z)$ and some applications}, Int.\ Math.\ Res.\ Notices 2009, no.\ 8, 1488--1545
\bibitem[S3]{S3} N.\ R.\ Scheithauer, {\em Some constructions of modular forms for the Weil representation of $\SL_2(\Z)$}, Nagoya Math.\ J.\ {\bf 220} (2015), 1--43
\bibitem[S4]{S4} N.\ R.\ Scheithauer, {\em Automorphic products of singular weight}, Compositio Math.\ {\bf 153} (2017), 1855--1892
  % \bibitem[S]{S} B.\ Schoeneberg, {\em Das Verhalten von mehrfachen Thetareihen bei Modulsubstitutionen}, Math.\ Ann. {\bf 116} (1939), 511--523
\bibitem[Sk1]{Sk1} N.-P.\ Skoruppa, {\em On the connection between Jacobi forms and modular forms of half-integral weight}, Dissertation,
  % Rheinische Friedrich-Wilhelms-Universit\"at, Bonn, 1984.
  Bonn Mathematical Publications, {\bf 159}, Bonn, 1985
\bibitem[Sk2]{Sk2} N.-P.\ Skoruppa, {\em Jacobi forms of critical weight and Weil representations}, Modular forms on Schiermonnikoog, 239--266, Cambridge Univ.\ Press, Cambridge, 2008
\bibitem[W]{W} A.\ Weil, {\em Sur certains groupes d'op\'{e}rateurs unitaires}, Acta Math. {\bf 111} (1964), 143--211
\bibitem[Wr]{Wr} F.\ Werner, {\em Vector valued Hecke theory}, Ph.D.\ Thesis, 2014, Tech\-ni\-sche Universit\"at Darmstadt, available at
https:/$\!$/tuprints.ulb.tu-darmstadt.de/4238/
\bibitem[Z]{Z} S.\ Zemel, {\em Integral bases and invariant vectors for {W}eil representations}, Res.\ Number Theory {\bf 9} (2023), Article number 5, 27 pp.\
%\bibitem[Z]{Z} W.\ Zhang, {\em Weil representation and arithmetic fundamental lemma}, Ann.\ of Math.\ (2) {\bf 193} (2021), 863--978 
\end{thebibliography}
\end{document}